\documentclass[11pt]{amsart}
\usepackage{amsopn}
\usepackage{amssymb, amscd}
\usepackage{multirow}
\usepackage{graphicx, graphics, epsfig}
\usepackage{faktor}
\usepackage{enumerate}
\usepackage{tikz}
\usepackage{pgfplots}
\usepackage{hyperref}
\usepackage{color}
\usepackage{pdfpages}

\newcommand\black{\color{black}}

\newcommand{\nc}{\newcommand}

\nc{\fg}{\mathfrak{f} } \nc{\vg}{\mathfrak{v} } \nc{\wg}{\mathfrak{w} }
\nc{\zg}{\mathfrak{z} } \nc{\ngo}{\mathfrak{n} } \nc{\kg}{\mathfrak{k} }
\nc{\mg}{\mathfrak{m} } \nc{\bg}{\mathfrak{b} } \nc{\ggo}{\mathfrak{g} } \nc{\eg}{\mathfrak{e} }
\nc{\ggob}{\overline{\mathfrak{g}} } \nc{\sog}{\mathfrak{so} }
\nc{\sug}{\mathfrak{su} } \nc{\spg}{\mathfrak{sp} } \nc{\slg}{\mathfrak{sl} }
\nc{\glg}{\mathfrak{gl} } \nc{\cg}{\mathfrak{c} } \nc{\rg}{\mathfrak{r} }
\nc{\hg}{\mathfrak{h} } \nc{\tg}{\mathfrak{t} } \nc{\ug}{\mathfrak{u} }
\nc{\dg}{\mathfrak{d} } \nc{\ag}{\mathfrak{a} } \nc{\pg}{\mathfrak{p} }
\nc{\sg}{\mathfrak{s} } \nc{\affg}{\mathfrak{aff} } \nc{\qg}{\mathfrak{q} } \nc{\lgo}{\mathfrak{l} }

\nc{\pca}{\mathcal{P}} \nc{\nca}{\mathcal{N}} \nc{\lca}{\mathcal{L}}
\nc{\oca}{\mathcal{O}} \nc{\mca}{\mathcal{M}} \nc{\tca}{\mathcal{T}}
\nc{\aca}{\mathcal{A}} \nc{\cca}{\mathcal{C}} \nc{\gca}{\mathcal{G}}
\nc{\sca}{\mathcal{S}} \nc{\hca}{\mathcal{H}} \nc{\bca}{\mathcal{B}}
\nc{\dca}{\mathcal{D}} \nc{\eca}{\mathcal{E}} \nc{\wca}{\mathcal{W}}

\nc{\vp}{\varphi} \nc{\ddt}{\tfrac{d}{dt}} \nc{\dsdt}{\tfrac{d^2}{dt^2}} \nc{\dds}{\frac{d}{ds}}
\nc{\dpar}{\frac{\partial}{\partial t}} \nc{\im}{\mathrm{i}}

\nc{\SO}{\mathrm{SO}} \nc{\Spe}{\mathrm{Sp}} \nc{\Sl}{\mathrm{SL}}
\nc{\SU}{\mathrm{SU}} \nc{\Or}{\mathrm{O}} \nc{\U}{\mathrm{U}} \nc{\Gl}{\mathrm{GL}}
\nc{\Se}{\mathrm{S}} \nc{\Cl}{\mathrm{Cl}} \nc{\Spin}{\mathrm{Spin}}
\nc{\Pin}{\mathrm{Pin}} \nc{\G}{\mathrm{GL}_n(\RR)} \nc{\g}{\mathfrak{gl}_n(\RR)}

\nc{\RR}{{\Bbb R}} \nc{\HH}{{\Bbb H}} \nc{\CC}{{\Bbb C}} \nc{\ZZ}{{\Bbb Z}}
\nc{\FF}{{\Bbb F}} \nc{\NN}{{\Bbb N}} \nc{\QQ}{{\Bbb Q}} \nc{\PP}{{\Bbb P}} \nc{\OO}{{\Bbb O}}

\nc{\vs}{\vspace{.2cm}} \nc{\vsp}{\vspace{1cm}} \nc{\ip}{\langle\cdot,\cdot\rangle}
\nc{\ipp}{(\cdot,\cdot)} \nc{\la}{\langle} \nc{\ra}{\rangle} \nc{\unm}{\tfrac{1}{2}}
\nc{\unc}{\tfrac{1}{4}} \nc{\und}{\tfrac{1}{16}} \nc{\no}{\vs\noindent}
\nc{\lam}{\Lambda^2(\RR^n)^*\otimes\RR^n} \nc{\tangz}{{\rm T}^{\rm Zar}}
\nc{\nor}{{\sf n}}  \nc{\mum}{/\!\!/} \nc{\kir}{/\!\!/\!\!/}
\nc{\Ri}{\tfrac{4\Ric_{\mu}}{||\mu||^2}} \nc{\ds}{\displaystyle}
\nc{\ben}{\begin{enumerate}} \nc{\een}{\end{enumerate}} \nc{\f}{\frac}
\nc{\lb}{[\cdot,\cdot]} \nc{\isn}{\tfrac{1}{||v||^2}}
\nc{\gkp}{(\ggo=\kg\oplus\pg,\ip)} \nc{\ukh}{(\ug=\kg\oplus\hg,\ip)}
\nc{\tgkp}{(\tilde{\ggo}=\kg\oplus\pg,\ip)}
\nc{\wt}{\widetilde}
\nc{\iop}{\mathtt{i}} \nc{\jop}{\mathtt{j}}

\nc{\Hess}{\operatorname{Hess}} \nc{\ad}{\operatorname{ad}}
\nc{\Ad}{\operatorname{Ad}} \nc{\rank}{\operatorname{rk}}
\nc{\Irr}{\operatorname{Irr}} \nc{\End}{\operatorname{End}}
\nc{\Aut}{\operatorname{Aut}} \nc{\Inn}{\operatorname{Inn}}
\nc{\Der}{\operatorname{Der}} \nc{\Ker}{\operatorname{Ker}}
\nc{\Iso}{\operatorname{Iso}} \nc{\Diff}{\operatorname{Diff}}
\nc{\Lie}{\operatorname{L}} \nc{\tr}{\operatorname{tr}} \nc{\dif}{\operatorname{d}}
\nc{\sen}{\operatorname{sen}} \nc{\modu}{\operatorname{mod}}
\nc{\CRic}{\operatorname{PP}} \nc{\Cric}{\operatorname{P}} \nc{\Ricci}{\operatorname{Ric}}
\nc{\sym}{\operatorname{sym}} \nc{\herm}{\operatorname{herm}} \nc{\symac}{\operatorname{sym^{ac}}}
\nc{\symc}{\operatorname{sym^{c}}} \nc{\scalar}{\operatorname{Sc}}
\nc{\grad}{\operatorname{grad}} \nc{\ricci}{\operatorname{Rc}} \nc{\kil}{\operatorname{B}} \nc{\cas}{\operatorname{C}} \nc{\lic}{\operatorname{L}}
\nc{\Nor}{\operatorname{Norm}}  \nc{\ricc}{\operatorname{Rc^{c}}}
\nc{\Ricc}{\operatorname{Ric^{c}}} \nc{\ricac}{\operatorname{Rc^{ac}}}
\nc{\Ricac}{\operatorname{Ric^{ac}}} \nc{\Riem}{\operatorname{Rm}} \nc{\Sec}{\operatorname{Sec}}
\nc{\riccig}{\operatorname{ric^{\gamma}}} \nc{\mm}{\operatorname{m}} \nc{\Mm}{\operatorname{M}}
\nc{\Le}{\operatorname{L}} \nc{\tang}{\operatorname{T}}
\nc{\level}{\operatorname{level}} \nc{\rad}{\operatorname{r}}
\nc{\abel}{\operatorname{ab}} \nc{\CH}{\operatorname{CH}} \nc{\Cone}{{\mathcal C}} \nc{\CCone}{\operatorname{CC}} \nc{\CP}{{\mathcal P}}
\nc{\mcc}{\operatorname{mcc}} \nc{\Adj}{\operatorname{Adj}}
\nc{\Order}{\operatorname{O}}  \nc{\inj}{\operatorname{inj}} \nc{\proy}{\operatorname{pr}}
\nc{\vol}{\operatorname{vol}} \nc{\Diag}{\operatorname{Dg}} \nc{\Diagg}{\operatorname{Diag}}
\nc{\Spec}{\operatorname{Spec}} \nc{\Ima}{\operatorname{Im}} \nc{\Rea}{\operatorname{Re}}
\nc{\spann}{\operatorname{span}} \nc{\Aff}{\operatorname{Aff}} \nc{\E}{\operatorname{E}} \nc{\id}{\operatorname{id}} \nc{\dete}{\operatorname{det}} \nc{\Crit}{\operatorname{Crit}} \nc{\val}{\operatorname{val}}

\theoremstyle{plain}
\newtheorem{theorem}{Theorem}[section]

\theoremstyle{definition}
\newtheorem{definition}[theorem]{Definition}

\theoremstyle{remark}
\newtheorem{remark}[theorem]{Remark}
\newtheorem{example}[theorem]{Example}

\title{On the stability of homogeneous Einstein manifolds II}

\author{Jorge Lauret}  \author{Cynthia Will}

\address{FaMAF, Universidad Nacional de C\'ordoba and CIEM, CONICET (Argentina)}
\email{jorgelauret@unc.edu.ar}  \email{cynthia.will@unc.edu.ar}

\thanks{This research was partially supported by a grant from Univ. Nac. de C\'ordoba, Argentina.}

\date{\today}

\begin{document}

\maketitle

\begin{abstract}
For any $G$-invariant metric on a compact homogeneous space $M=G/K$, we give a formula for the Lichnerowicz Laplacian restricted to the space of all $G$-invariant symmetric $2$-tensors in terms of the structural constants of $G/K$.  As an application, we compute the $G$-invariant spectrum of the Lichnerowicz Laplacian for all the Einstein metrics on most generalized Wallach spaces and any flag manifold with $b_2(M)=1$.  This allows to deduce the $G$-stability and critical point types of each of such Einstein metrics as a critical point of the scalar curvature functional.  
\end{abstract}

\tableofcontents


\section{Introduction}\label{intro}

Back in 1915, Hilbert considered in \cite{Hlb} the simplest curvature functional on the space $\mca_1$ of all unit volume Riemannian metrics on a given compact manifold $M$, i.e., the total scalar curvature 
\begin{equation}\label{sct}
\widetilde{\scalar}(g):=\int_M \scalar(g)\; d\vol_g,
\end{equation}
and proved that their critical points are precisely the metrics of constant Ricci curvature, i.e., $\ricci(g)=\rho g$ for some $\rho\in\RR$, so-called {\it Einstein} metrics (see \cite[4.21]{Bss}).  A long time later, in the seventies, Berger \cite{Brg} and Koiso \cite{Kso2} independently discovered that if one considers the restriction of $\widetilde{\scalar}$ to the space $\cca_1$ of all unit volume constant scalar curvature metrics, then the coindex and nullity of any Einstein metric are both finite (modulo trivial variations).  A fundamental problem is therefore to determine whether a given Einstein metric $g$ is {\it stable} (or linearly stable), i.e., coindex and nullity indeed vanish, or equivalently, the Hessian of $\widetilde{\scalar}$ is negative definite on the tangent space $T_g\cca_1$ (modulo trivial variations).  Any stable Einstein metric is in particular a local maximum of $\widetilde{\scalar}|_{\cca_1}$ and it is {\it rigid}, in the sense that every Einstein metric sufficiently close to $g$ is isometric to $g$ up to scaling.  It was proved by Koiso that the non-degeneracy of the Hessian is actually enough to get rigidity (see \cite[Proposition 3.3]{Kso} and \cite[12.66]{Bss}).

Stability seems to be a very strong property, which is in addition hard to be established.  Indeed, the only known examples so far of stable Einstein metrics (or local maxima of $\widetilde{\scalar}|_{\cca_1}$) with positive scalar curvature are given by some irreducible symmetric spaces.  The space $T_g\cca_1$ coincides, modulo trivial variations, with $\tca\tca_g=\Ker\delta_g\cap\Ker\tr_g$, the space of divergence-free (or transversal) and traceless symmetric $2$-tensors, so-called {\it TT-tensors}, and if $\ricci(g)=\rho g$, then for any $T\in\tca\tca_g$ the Hessian is given by
\begin{equation}\label{LL-intro}
\widetilde{\scalar}''_g(T,T) = \unm\la(2\rho\id-\Delta_L)T,T\ra_g,
\end{equation}
where $\Delta_L$ is the {\it Lichnerowicz Laplacian} of $g$ (see \cite[4.64]{Bss}).  Thus $g$ is stable if and only if $2\rho<\lambda_L$, where $\lambda_L$ denotes the smallest eigenvalue of $\Delta_L|_{\tca\tca_g}$ (cf.\ \cite[\S 4]{CaoHe} and \cite[\S 1]{WngWng}).  The number $\lambda_L$ is usually hard to compute or even estimate, it is known for only few metrics.  Note that $\widetilde{\scalar}''_g|_{\tca\tca_g}$ is non-degenerate if and only if $2\rho\notin\Spec\left(\Delta_L|_{\tca\tca_g}\right)$.   

We assume that $M$ is homogeneous in this paper.  After fixing a transitive action of a compact Lie group $G$ on $M$, we study the stability of $G$-invariant Einstein metrics, which are precisely the critical points of 
$
\scalar:\mca_1^G\longrightarrow \RR,  
$ 
where $\mca^G_1$ is the finite-dimensional manifold of all unit volume $G$-invariant metrics on $M$.  The $G$-action provides a presentation $M=G/K$ of $M$ as a homogeneous space, where $K\subset G$ is the isotropy subgroup at some origin point $o\in M$.  It is proved in \cite[\S 3.4]{stab-tres} that 
$$
T_g\mca_1^G = T_gN_G(K)^*g \oplus \tca\tca_g^G, 
$$
where $N_G(K)\subset\Diff(M)$ is the normalizer of $K$ and $\tca\tca_g^G$ is the space of $G$-invariant TT-tensors.  An Einstein metric $g\in\mca_1^G$ is therefore called  $G$-{\it stable} when $\scalar''_g|_{\tca\tca_g^G}<0$, or equivalently, $2\rho <\lambda_L^G$, where $\lambda_L^G$ is the smallest eigenvalue of $\Delta_L|_{\tca\tca_g^G}$.  Note that 
$$
\lambda_L\leq\lambda_L^G. 
$$  
One of the aims of \cite{stab-tres} and this paper is to show that this `algebraic' upper bound for $\lambda_L$ is in most cases enough to establish the instability of a homogeneous Einstein metric.      
  
Any $G$-stable Einstein metric is in particular a local maximum of $\scalar:\mca_1^G\longrightarrow \RR$, and if $\scalar''_g|_{\tca\tca_g^G}$ is non-degenerate, then $g$ is $G$-{\it rigid}, in the sense that $g$ is isolated in the moduli space of $G$-invariant Einstein metrics on $M$ up to pull-back by $N_G(K)$ and scaling.  $G$-stability is of course a necessary condition for the classical stability described above, and it is also extremely rare.  As far as we know, besides irreducible symmetric metrics and the special case when any subalgebra of $\ggo$ containing $\kg$ is of the form $\kg\oplus\ag$ with $[\ag,\ag]=0$ (e.g., if $K$ is a maximal subgroup of $G$, see \cite{WngZll, Bhm}), the only known $G$-stable Einstein metrics with $\dim{\mca_1^G}>1$ are:
\begin{enumerate}[{\small $\bullet$}] 
\item the standard metric on $\SU(2)$, $E_7/\SO(8)$ and $E_8/\Spin(8)\times\Spin(8)$, discovered in \cite{AbvArvNknSss} (see also \cite[Remark 2.3]{ChnNkn}), 

\item and the unique K\"ahler-Einstein metric on each of the $13$ flag manifolds with second Betti number $b_2(M)$ equal to $1$, which can be deduced from \cite[Theorem 3.1]{AnsChr} (here $G$ is always an exceptional simple Lie group).     
\end{enumerate}
It is unknown whether these $G$-stable Einstein metrics are stable, as well as whether they realize the Yamabe invariant of $M$ (i.e., whether $\scalar(g)$ is the supreme among all Yamabe metrics of $M$, which are those with the smallest scalar curvature in its unit volume conformal class).  

\begin{remark}\label{new}
After the first version of the present paper was uploaded to arXiv, many new examples of $G$-stable standard metrics have been found in \cite{stab} and the $G$-stable Einstein metric on $E_7/\SO(8)$ has been proved to be stable in \cite{SchSmmWng}.  
\end{remark}

On the other hand, if $g$ is $G$-{\it unstable}, i.e., $\lambda_L^G<2\rho$ (or $\scalar''_g(T,T)>0$ for some $T\in{\tca\tca_g^G}$), then $g$ is also dynamically unstable for the normalized Ricci flow, in the sense that there is an ancient solution emerging from $g$ (see \cite{CaoHe, Krn}), and $g$ does not realize the Yamabe invariant of $M$ (see \cite[Theorem 5.1]{BhmWngZll}).  Excepting the ones described in the two items above, all the homogeneous Einstein metrics studied in this paper were found to be $G$-unstable.  Indeed, $G$-instability is considered a generic property of a compact homogeneous Einstein manifold by the experts.   

It follows from \eqref{LL-intro} that the spectrum of the Lichnerowicz Laplacian $\Delta_L$ restricted to $\tca\tca_g^G$ plays a crucial role in the study of $G$-stability and $G$-non-degeneracy.  Consider any reductive decomposition $\ggo=\kg\oplus\pg$ of $M=G/K$ and all the corresponding usual identifications $T_oM\equiv \pg$, $\sca^2(M)^G\equiv \sym(\pg)^K$, etc.  According to \cite[Lemma 4.7]{stab-tres}, for any $g\in\mca^G$,  
$$
\Delta_LT = \la\lic_\pg A\cdot,\cdot\ra, \qquad\forall T\in\tca\tca_g^G, \quad T=\la A\cdot,\cdot\ra\in \sca^2(M)^G, \quad A\in\sym(\pg)^K,
$$
where
$
\lic_\pg=\lic_\pg(g):\sym(\pg)\longrightarrow\sym(\pg) 
$
is the self-adjoint operator defined by
\begin{equation}\label{Lp-intro}
\la \lic_\pg A,B\ra = \unm \la\theta(A)\mu_\pg,\theta(B)\mu_\pg\ra + 2\tr{\Mm_{\mu_\pg}AB}, \qquad\forall A,B\in\sym(\pg),
\end{equation}
and $\ip$ is determined by $g$.  Here $\mu_\pg:= \proy_\pg\circ \mu|_{\pg\times\pg} :\pg\times\pg\longrightarrow\pg$ and $\Mm:\Lambda^2\pg^*\otimes\pg\rightarrow\sym(\pg)$ is the moment map for the usual representation $\theta$ of $\glg(\pg)$ on $\Lambda^2\pg^*\otimes\pg$ (see \S\ref{ric-sec}).  Thus $\lambda_L^G$ (denoted by $\lambda_\pg$ in subsequent sections) is the smallest eigenvalue of the operator $\lic_\pg$ restricted to $\sym(\pg)^K$ (modulo trivial variations).  This was applied to the naturally reductive case in \cite{stab-tres}, where formula \eqref{Lp-intro} considerably simplifies.    

In this paper, we give a formula for the Lichnerowicz Laplacian $\lic_\pg$ in the general case, in terms of the well-known structural constants $[ijk]$'s of $G/K$ relative to a $Q$-orthogonal decomposition 
$\pg=\pg_1\oplus\dots\oplus\pg_r$ in $\Ad(K)$-invariant subspaces, where $Q$ is any bi-invariant inner product on $\ggo$ (see \eqref{ijk}).  Assume that the subspaces $\pg_k$'s are all $\Ad(K)$-irreducible and pairwise inequivalent, which implies that   
$\left\{ \tfrac{1}{\sqrt{d_1}}I_1,\dots, \tfrac{1}{\sqrt{d_r}}I_r\right\}$
is an orthonormal basis of $\sym(\pg)^K$, where $d_k:=\dim{\pg_k}$ and $I_k|_{\pg_i}:=\delta_{ki}I$.

\begin{theorem}\label{formLp-intro} (See Theorem \ref{formLp}).   
For any metric $g=x_1Q|_{\pg_1}+\dots+x_rQ|_{\pg_r}\in\mca^G$, the entries of the matrix of $\lic_\pg$ with respect to the above basis are given by
$$
[\lic_\pg]_{kk} = \tfrac{1}{d_k}\sum_{i,j\ne k} \tfrac{x_k}{x_ix_j}[ijk]  + \tfrac{1}{d_k}\sum_{i\ne k}  \tfrac{x_i}{x_k^2}[ikk], \quad
[\lic_\pg]_{km} = \tfrac{1}{\sqrt{d_k}\sqrt{d_m}} \sum_{i} \tfrac{x_i^2-x_k^2-x_m^2}{x_ix_kx_m} [ikm].  
$$
\end{theorem}

The formula is also useful beyond the multiplicity-free case, see Remarks \ref{rem-Lp} and \ref{Lp-gen}.  Theorem \ref{formLp-intro} allows the systematic computation of the spectrum of $\lic_\pg$ for any homogeneous Einstein metric available in the literature, provided that the numbers $d_k$'s and $[ijk]$'s are known for the homogeneous space.  This is usually the case since Einstein equations are often written in terms of these numbers (see \S\ref{E-sec}).    

The rest of the paper is explorative in nature.  In \S\ref{W-sec}, we use the formula in Theorem \ref{formLp-intro} to compute the spectrum of the Lichnerowicz Laplacian $\lic_\pg$ for all Einstein metrics on any generalized Wallach space $G/K$ with $G$ simple ($5$ infinite families and $10$ exceptional examples; each space admits between $2$ and $4$ Einstein metrics), except $\SO(k+l+m)/\SO(k)\times\SO(l)\times\SO(m)$ and $\Spe(k+l+m)/\Spe(k)\times\Spe(l)\times\Spe(m)$, where $k,l,m$ are pairwise different.  We also apply Theorem \ref{formLp-intro} in \S\ref{flag-sec} to do the same for all Einstein metrics on flag manifolds with $b_2(M)=1$ ($13$ exceptional examples, each one admitting between $3$ and $6$ Einstein metrics).  The $G$-stability and critical point types of each of these metrics so obtained are collected in Tables \ref{W2}, \ref{W4}, \ref{W5} and Tables \ref{FS3}-\ref{FS6}, respectively.  The most noteworthy finds of our exploration are:     

\begin{enumerate}[{\small $\bullet$}] 
\item The only local maxima of $\scalar|_{\mca_1^G}$ are the standard metrics on $\SU(2)$, $E_7/\SO(8)$ and $E_8/\Spin(8)\times\Spin(8)$, and the K\"ahler metrics on the flag manifolds with $b_2(M)=1$.  However, the K\"ahler-Einstein metrics on flag manifolds with $b_2(M)=2$ and three isotropy summands are never local maxima (cf.\ \cite{CaoHI}).    

\item Any standard Einstein metric on a generalized Wallach space ($3$ infinite families and $7$ isolated examples) is either a local maximum or a local minimum of $\scalar|_{\mca_1^G}$.  The value of the volume normalized scalar curvature $\scalar_N$ at the standard metric is always strictly less than for all the other Einstein metrics, except precisely in the cases when it is a local maximum.    

\item A generalized Wallach space admits a local minimum if and only if it admits four Einstein metrics.  

\item For any of the three K\"ahler-Einstein metrics on any of the generalized Wallach spaces which are also flag manifolds (necessarily with $b_2(M)=2$), the kernel of the Lichnerowicz Laplacian $\Delta_L|_{\tca\tca_g^G}$ is non-trivial.  By \cite{PRP}, this implies that there exists a symmetric $2$-tensor $T$ as close to $\ricci(g)$ as you like such that either the existence or the (local) uniqueness of a solution to the prescribed Ricci curvature problem $\ricci(g')=cT$ fail (see \S\ref{RLI}).  Indeed, for each of these K\"ahler-Einstein metrics $g_{KE}$, we provide a pairwise non-homothetic curve $g_t$, $t\in [1,\infty)$ of $G$-invariant metrics ($\scalar(g_t)$ is strictly increasing) having the same volume as $g_1=g_{KE}$ and such that $\ricci(g_t)=\ricci(g_{KE})$ for all $t\geq 1$.  This only happens for these metrics, all the other ones are Ricci locally invertible.   

\item All the Einstein metrics studied are non-degenerate critical points of $\scalar|_{\mca_1^G}$ and in particular $G$-rigid.  The $G$-rigidity was previously known, it follows from the finiteness (up to scaling) of the set of Einstein metrics on each of these homogeneous spaces (see \cite{FrsLmsNkn, Arv}).  

\item Any non-K\"ahler Einstein metric on a flag manifold with $b_2(M)=1$ is a saddle point of coindex $1$, except for one of the metrics on each of the three spaces
$E_8/\SU(10)\times\SU(3)\times\SU(1)$, $E_8/ \SU(5)\times\SU(4)\times\U(1)$ and $E_8/ \SU(5)\times \SU(3)\times \SU(2)\times\U(1)$, which has coindex $2$.  
\end{enumerate}

The normalized Ricci flow on $\mca_1^G$ is precisely the gradient flow of $\scalar$, whose dynamical behavior is mostly governed by the $G$-stability types of their fixed points, the $G$-invariant Einstein metrics.  The Ricci flow on generalized Wallach spaces was studied in \cite{AbvArvNknSss}, where the authors prove that, generically, there are only three possible dynamical behaviors.  Our approach provides an alternative proof for that.  On flag manifolds with $b_2(M)=1$, the dimension of the unstable manifold for each fixed point at the infinity is given in \cite[Theorem 3.1]{AnsChr}.  It can be shown that such dimension is precisely one more than the coindex of the corresponding critical point.  Our computations of the coindex agree with these results.      
 
\begin{remark} (cf.\ \cite[Remark 1.6]{CaoHe}).  
As proposed in \cite{GbbHrtPp}, higher-dimensional versions of the spherically symmetric Schwarzschild black hole can be constructed by replacing the $2$-sphere with other compact Einstein manifolds.  The instability condition for the $(n+2)$-dimensional black hole obtained is given by 
$$
\lambda_L < \frac{\rho}{n-1}\left(4-\frac{(n-5)^2}{4}\right) = \rho\left(\frac{9-n}{4}\right), 
$$  
where $\lambda_L$ and $\rho$ correspond to the Einstein manifold $M^n$.  It is also proved in \cite{GbbHrtPp} that this instability criterion is identical to that for the instability of the Freund-Rubin compactification $AdS_m\times M^n$ in the context of the AdS/CFT correspondence.  Remarkably, among all the homogeneous Einstein manifolds studied in this paper, in spite of $\lambda_L^G$ is in many cases negative (recall that $\lambda_L\leq\lambda_L^G$), the Einstein-K\"ahler metric on the full flag manifold $M^6=\SU(3)/\Se(\U(1)^3)$ (which is not the standard metric) is the only one satisfying the instability condition above with $\lambda_L$ replaced by $\lambda_L^G$.  Since $\lambda_L^G=0$ for this particular metric (see Table \ref{W2}, line $W_2$, $k=1$), it follows that $\lambda_L\leq 0$, providing an $8$-dimensional generalized black hole which is unstable.  Topologically, $M^6$ is the projective holomorphic tangent bundle of $\CC P^2$ (see \cite[Example 13.80 and Chapter 8]{Bss}).          
\end{remark} 
  
\vs \noindent {\it Acknowledgements.}  We are very grateful with Ioannis Chrysikos, Lino Grama, Yurii Nikonorov and Wolfgang Ziller for many helpful conversations.

\section{Preliminaries}\label{preli}

Let $M$ be a compact connected differentiable manifold.  We assume that $M$ is homogeneous and fix an almost-effective transitive action of a compact Lie group $G$ on $M$.  The $G$-action determines a presentation $M=G/K$ of $M$ as a homogeneous space, where $K\subset G$ is the isotropy subgroup at some point $o\in M$.  Let $\mca^G$ denote the finite-dimensional manifold of all $G$-invariant Riemannian metrics on $M$.   

It is well known that $g\in\mca_1^G$ is Einstein, i.e., $\ricci(g)=\rho g$ for some $\rho\in\RR$ (which is necessarily positive if $G$ is non-abelian), if and only if $g$ is a critical point of the scalar curvature functional $\scalar:\mca_1^G\longrightarrow \RR$, where $\mca^G_1\subset\mca^G$ is the codimension one submanifold of all unit volume metrics.  We study in this paper the stability of $G$-invariant Einstein metrics on $M$ as critical points of $\scalar|_{\mca_1^G}$, which is encoded in the signature of the second derivative or Hessian $\scalar''$.  We refer to \cite{stab-tres} for a more detailed treatment.

\subsection{$G$-stability}\label{Gstab-sec}
The tangent space of $\mca_1^G$ at $g$ is given by the subspace of $\sca^2(M)^G$ of $g$-traceless elements, where $\sca^2(M)^G$ is the space of all $G$-invariant symmetric $2$-tensors.  It follows from \cite[Corollary 3.12]{stab-tres} that 
\begin{equation}\label{tang}
T_g\mca_1^G = T_g\Aut(G/K)\cdot g \oplus \tca\tca_g^G, 
\end{equation}
where $\Aut(G/K)\subset\Diff(M)$ is the Lie group of automorphisms of $G$ taking $K$ onto $K$, acting by pullback on $\mca^G$ (i.e., $T_g\Aut(G/K)\cdot g$ is the space of trivial variations of $g$), and $\tca\tca_g^G:=(\Ker\delta_g\cap\Ker\tr_g)^G$ is the space of $G$-invariant TT-tensors.  Note that since $G$ is compact, $T_g\Aut(G/K)\cdot g=T_gN_G(K)\cdot g$, where $N_G(K)$ is the normalizer of $K$ in $G$ and so $T_gN_G(K)\cdot g$ is the tangent space of the $G$-equivariant isometry class of $g$ (see \cite[\S3.3]{stab-tres}).   

\begin{definition}\label{stab-def}
An Einstein metric $g\in\mca^G_1$ is said to be,  
\begin{enumerate}[{\small $\bullet$}] 
\item {\it $G$-stable}: $\scalar''_g|_{\tca\tca_g^G}<0$ (in particular, $g$ is a local maximum of $\scalar|_{\mca^G_1}$). 

\item {\it $G$-unstable}: $\scalar''_g(T,T)>0$ for some $T\in \tca\tca_g^G$ ($g$ is a saddle point, unless $\scalar''_g|_{\tca\tca_g^G}>0$, in which case $g$ is a local minimum of $\scalar|_{\mca^G_1}$).  The {\it coindex} is the dimension of the maximal subspace of $\tca\tca_g^G$ on which $\scalar''_g$ is positive definite.   

\item {\it $G$-non-degenerate}: $\scalar''_g|_{\tca\tca_g^G}$ is non-degenerate (thus $g$ is $G$-rigid, in the sense that it is an isolated critical point up to the $\Aut(G/K)$-action), and otherwise, {\it $G$-degenerate}.   
\end{enumerate}
\end{definition}

Recall that without assuming $G$-invariance, the classical stability of $g$ means that $\scalar''_g$ is negative definite on $\tca\tca_g$, the infinite dimensional vector space of all unit volume constant scalar curvature (non-trivial) variations of $g$. 

It is worth noticing that the space of trivial variations vanishes, i.e., $T_g\mca_1^G = \tca\tca_g^G$, under any of the following conditions: $g$ is naturally reductive with respect to $G$; the isotropy representation does not contain the trivial representation; $G$ is semisimple and the isotropy representation is {\it mutiplicity-free} (i.e., any two different $\Ad(K)$-invariant irreducible subspaces are inequivalent as $\Ad(K)$-representations), e.g., when $\rank(G)=\rank(K)$ (see \cite[\S3.3]{stab-tres}).

\subsection{Ricci curvature}\label{ric-sec}
Given a compact homogeneous space $M=G/K$ as above, we consider any reductive decomposition $\ggo=\kg\oplus\pg$, giving rise to the usual identification $T_oM\equiv\pg$, where $\ggo$ and $\kg$ are the Lie algebras of $G$ and $K$, respectively.  Let $\mu$ denote the Lie bracket of $\ggo$.  We extend $\ip:=g_o$ in the usual way to inner products on $\glg(\pg)$ and $\Lambda^2\pg^*\otimes\pg$, respectively:
$$
\la A,B\ra:= \tr{AB^t}, \qquad \la\lambda,\lambda\ra:=\sum |\ad_\lambda{X_i}|^2 = \sum |\lambda(X_i,X_j)|^2,
$$
where $\{ X_i\}$ is any orthonormal basis of $\pg$ relative to $\ip$.  If
\begin{equation}\label{mup}
\mu_\pg:= \proy_\pg\circ \mu|_{\pg\times\pg} :\pg\times\pg\longrightarrow\pg,
\end{equation}
where $\proy_\pg:\ggo\rightarrow\pg$ is the projection on $\pg$ relative to $\ggo=\kg\oplus\pg$, then  the Ricci operator $\Ricci(g)$ of the metric $g$ is given by
\begin{equation}\label{Ric2}
\Ricci(g) = \Mm_{\mu_\pg} - \unm\kil_\mu,
\end{equation}
where $\la\kil_\mu\cdot,\cdot\ra:=\kil_{\ggo}|_{\pg\times\pg}$, $\kil_{\ggo}$ denotes the Killing form of the Lie algebra $\ggo$ and
\begin{equation}\label{mm2}
\la\Mm_{\mu_\pg},A\ra := \unc\la\theta(A)\mu_\pg,\mu_\pg\ra, \qquad\forall A\in\glg(\pg).
\end{equation}
Here $\theta$ is the representation of $\glg(\pg)$ given by,
\begin{equation}\label{tita}
\theta(A)\lambda := A\lambda(\cdot,\cdot) - \lambda(A\cdot,\cdot) - \lambda(\cdot,A\cdot), \qquad \forall A\in\glg(\pg), \quad \lambda\in\Lambda^2\pg^*\otimes\pg.
\end{equation}
The function $\Mm:\Lambda^2\pg^*\otimes\pg\rightarrow\sym(\pg)$ is therefore the {\it moment map} from geometric invariant theory (see e.g.\ \cite{BhmLfn} and the references therein) for the representation $\theta$ of $\glg(\pg)$.  Alternatively, the moment map can be defined by,
\begin{equation}\label{mm4}
\la\Mm_{\mu_\pg}X,X\ra = -\unm\sum \la\mu_\pg(X,X_i),X_j\ra^2+ \unc\sum \la\mu_\pg(X_i,X_j),X\ra^2, \qquad\forall X\in\pg.
\end{equation} 
It follows from \eqref{mm2} and \eqref{tita} that
\begin{equation}\label{scal2}
\scalar(g) = -\unc|\mu_\pg|^2  - \unm\tr{\kil_\mu}.
\end{equation}
We refer to \cite{alek} for more details on this formula for the Ricci curvature from the moving bracket approach point of view.

\subsection{Lichnerowicz Laplacian}\label{Lp-sec}
It is well known that the second variation of the total scalar curvature at any Einstein metric $g$ on $M$, say with $\ricci(g)=\rho g$, is given on $\tca\tca_g$ by
\begin{equation}\label{ScLL}
\scalar''_g = \unm\la (2\rho\id-\Delta_L)\cdot,\cdot\ra_g,
\end{equation}
where $\Delta_L$ is the Lichnerowicz Laplacian of $g$ (see \cite[4.64]{Bss}).  Consider the self-adjoint operator
$$
\lic_\pg=\lic_\pg(g):\sym(\pg)\longrightarrow\sym(\pg), \qquad 
\sym(\pg):=\{A:\pg\rightarrow\pg:A^t=A\},
$$
defined by
\begin{equation}\label{Lp}
\la \lic_\pg A,B\ra = \unm \la\theta(A)\mu_\pg,\theta(B)\mu_\pg\ra + 2\tr{\Mm_{\mu_\pg}AB}.\qquad\forall A,B\in\sym(\pg).
\end{equation}
It is proved in \cite[Lemma 4.7]{stab-tres} that this is precisely the operator defined by $\Delta_L$ under the usual identifications, that is,  
$$
\Delta_LT = \la\lic_\pg A\cdot,\cdot\ra, \qquad\forall T\in\tca\tca_g^G, \quad T=\la A\cdot,\cdot\ra\in \sca^2(M)^G, \quad A\in\sym(\pg)^K,
$$
where 
$$
\sym(\pg)^K:=\{A\in\sym(\pg): [\Ad(K),A]=0\} \equiv \sca^2(M)^G.
$$
According to \eqref{ScLL}, the $G$-stability type of $g$ is therefore determined by how is the constant $2\rho$ suited relative to the spectrum of the Lichnerowicz Laplacian $\lic_\pg$.  If $\lambda_\pg=\lambda_{\pg}(g)$ and $\lambda_\pg^{max}=\lambda_{\pg}^{max}(g)$ denote, respectively, the minimum and maximum eigenvalue of  $\lic_\pg=\lic_\pg(g)$ restricted to the subspace $\tca\tca_g^G$, then an Einstein metric $g\in\mca_1^G$, say $\ricci(g)=\rho g$, is

\begin{enumerate}[{\small $\bullet$}] 
\item $G$-stable if and only if $2\rho<\lambda_\pg$,

\item $G$-unstable if and only if $\lambda_\pg<2\rho$,

\item $G$-non-degenerate if and only if $2\rho\notin\Spec\left(\lic_\pg |_{\tca\tca_g^G}\right)$, 

\item a local minimum of $\scalar|_{\mca_1^G}$ if $\lambda_\pg^{max}<2\rho$.  
\end{enumerate}
We note that on trivial variations (see \eqref{tang}), $\lic_\pg|_{T_g\Aut(G/K)\cdot g}=2\rho\id$ (see \cite[(30)]{stab-tres}).  Recall that $\lambda_\pg$ was denoted by $\lambda_L^G$ in \S\ref{intro}.  

Given an intermediate subgroup $K\subset H\subset N_G(K)$, if $\sym(\pg)^H$ is the subspace of $\sym(\pg)^K$ of those maps which are $\Ad(H)$-invariant, then 
$$
\lic_\pg\sym(\pg)^H\subset\sym(\pg)^H,
$$ 
for any $g\in\mca^{G, H}:=\mca^G\cap\sym(\pg)^H$, the submanifold of $G$-invariant metrics which are in addition $\Ad(H)$-invariant (see \cite[\S4.4]{stab-tres}).  Note that since $\lambda_\pg\leq\Spec(\lic_\pg|_{\tca\tca_g^G\cap\sym(\pg)^H})$, if the minimum eigenvalue of $\lic_\pg|_{\tca\tca_g^G\cap\sym(\pg)^H}$ is $\leq 2\rho$, then $g$ is $G$-unstable, and that $g$ is $G$-degenerate if $2\rho\in\Spec(\lic_\pg|_{\tca\tca_g^G\cap\sym(\pg)^H})$.   

\begin{remark}\label{H}
It may be the case that $\mca^G=\Aut(G/K)\cdot\mca^{G, H}$, that is, any $G$-invariant metric is isometric to an $\Ad(H)$-invariant one.  A typical example of this behavior is the Stiefel manifold $\SO(k+2)/\SO(k)$, $k\geq 3$ (see \cite[\S 4]{Krr} and Examples \ref{Stiefel} and \ref{Stiefel3} below).  Since the function $\scalar$ therefore takes the same values on $\mca_1^{G,H}$ than on $\mca_1^G$, the $G$-stability and critical point types of an Einstein metric $g\in\mca_1^{G,H}$ as a critical point of $\scalar|_{\mca_1^{G,H}}$ and $\scalar|_{\mca_1^G}$ coincide, and so such types can be obtained from the spectrum of $\lic_\pg|_{\tca\tca_g^G\cap\sym(\pg)^H}$.  
\end{remark}

\subsection{Structural constants}\label{str-sec}
We fix a bi-invariant inner product $Q$ on $\ggo$ and consider the $Q$-orthogonal reductive decomposition $\ggo=\kg\oplus\pg$, and as a background metric, $Q|_{\pg}\in\mca^G$.  For each $Q$-orthogonal decomposition,
\begin{equation}\label{dec2}
\pg=\pg_1\oplus\dots\oplus\pg_r,
\end{equation}
in $\Ad(K)$-invariant  subspaces (not necessarily $\Ad(K)$-irreducible), we denote the corresponding structural constants by
\begin{equation}\label{ijk}
[ijk]:=\sum_{\alpha,\beta,\gamma} Q([e_{\alpha}^i,e_{\beta}^j], e_{\gamma}^k)^2,
\end{equation}
where $\{ e_\alpha^i\}$, $\{ e_{\beta}^j\}$ and $\{ e_{\gamma}^k\}$ are $Q$-orthonormal basis of $\pg_i$, $\pg_j$ and $\pg_k$, respectively.  Note that $[ijk]$ is invariant under any permutation of $ijk$ and it vanishes if and only if $[\pg_i,\pg_j]\perp\pg_k$.  These numbers have been strongly used in the literature since McKenzie Wang and Wolfgang Ziller introduced them in their pioneering article \cite{WngZll}.  

For each $g\in\mca^G$, there exists at least one decomposition as in \eqref{dec2} which is also $g$-orthogonal and
\begin{equation}\label{metric}
g=x_1Q|_{\pg_1}+\dots+x_rQ|_{\pg_r}, \qquad x_i>0.
\end{equation}
The metric $g$ will often be denoted by $(x_1,\dots,x_r)$.  We note that 
\begin{equation}\label{nmup}
|\mu_\pg|^2 = \sum_{i,j,k} \tfrac{x_k}{x_i x_j}[ijk],
\end{equation}
as $\left\{\frac{1}{\sqrt{x_i}}e_\alpha^i\right\}$ is a $g$-orthonormal basis of $\pg_i$ for all $i$.

\subsection{Einstein equations}\label{E-sec}
Recall from \eqref{Ric2} that the Ricci operator of $g$ is given by $\Ricci(g)=\Mm_{\mu_\pg}-\unm\kil_\mu$.  If we assume that, for each $k$, $\Mm_{\mu_\pg}|_{\pg_k}=m_kI_{\pg_k}$ for some $m_k\in\RR$ ($I_\sg$ denotes the identity map on $\sg$ for any vector space $\sg$), then it follows from \eqref{mm4} that
\begin{equation}\label{mkstr}
m_k=-\tfrac{1}{2d_k}\sum_{i,j} \tfrac{x_j}{x_ix_k}[ijk]  +\tfrac{1}{4d_k}\sum_{i,j} \tfrac{x_k}{x_ix_j}[ijk], \qquad\forall k=1,\dots,r,
\end{equation}
where $d_k:=\dim{\pg_k}$ and $g=(x_1,\dots,x_r)$.  Let us also assume that $-\kil_\ggo|_{\pg_k\times\pg_k}=b_kQ$, $b_k\in\RR$, which implies that $\kil_\mu|_{\pg_k}=-\tfrac{b_k}{x_k}I_{\pg_k}$.  Note that $b_k\geq 0$, where equality holds if and only if $\pg_k\subset\zg(\ggo)$, and that if $G$ is semisimple and $Q=-\kil_\ggo$, then $b_k=1$ for all $k$.

We therefore obtain that
$$
\Ricci(g)|_{\pg_k}=\rho_kI_{\pg_k}, \qquad \rho_k=\tfrac{b_k}{2x_k}+m_k,
$$
where each number $\rho_k$ can be written in the following ways:
\begin{align}
\rho_k  =& \tfrac{b_k}{2x_k} -\tfrac{1}{2d_k}\sum_{i,j} \tfrac{x_j}{x_ix_k}[ijk]  +\tfrac{1}{4d_k}\sum_{i,j} \tfrac{x_k}{x_ix_j}[ijk], \notag \\
=& \tfrac{b_k}{2x_k} -\tfrac{1}{4d_k}\sum_{i,j} \left(\tfrac{x_i}{x_jx_k} +\tfrac{x_j}{x_ix_k} - \tfrac{x_k}{x_ix_j}\right)[ijk], \label{rick}  \\
=& \tfrac{b_k}{2x_k} -\tfrac{1}{4d_k}\sum_{i,j} \tfrac{x_i^2+x_j^2-x_k^2}{x_ix_jx_k}[ijk].  \notag
\end{align}
The Einstein equations therefore become: $\Ricci(g)=\rho I$ if and only if
$$
\rho_k=\rho, \quad\forall k=1,\dots,r \quad\mbox{and}\quad \la\Ricci(g)\pg_i,\pg_j\ra=0, \quad\forall i\ne j.
$$
The two assumptions we have made above (namely, $\Mm_{\mu_\pg}|_{\pg_k}=m_kI_{\pg_k}$ and $-\kil_\ggo|_{\pg_k}=b_kQ|_{\pg_k}$) hold, for instance, in the case when the $\pg_k$'s are all $\Ad(K)$-irreducible, or more in general, if they are all $\Ad(H)$-invariant and irreducible for some intermediate subgroup $K\subset H\subset N_G(K)$.  If in addition the $\pg_k$'s are pairwise inequivalent as $K$-representations (resp.\ $H$-representations), then any $G$-invariant (resp.\ $(G\times H)$-invariant) metric is of the form \eqref{metric} and the right-hand side Einstein equations are automatically fulfilled.

These intermediate subgroups $K\subset H\subset N_G(K)$ can be used to search for Einstein metrics in the difficult case when $G/K$ is not multiplicity-free (see \cite{Stt} and references therein).     

\begin{example}
The isotropy representation $\pg$ of $\SO(7)/\SO(2)$ is highly non-multiplicity free, it decomposes in $10$ one-dimensional and $5$ two-dimensional irreducible $\SO(2)$-sub\-repre\-sen\-tations.  However, for the subgroup $H=\SO(4)\times\SO(2)$ of the normalizer of $\SO(2)$, $\pg$ only admits $4$ inequivalent $\Ad(H)$-irreducible subspaces of dimensions $2$, $4$, $6$ and $8$, respectively.   
\end{example}

\subsection{Scalar curvature}\label{Sc-sec}
It follows from \eqref{scal2} and \eqref{nmup} (or alternatively, from \eqref{rick}) that the scalar curvature of a metric $g=(x_1,\dots,x_r)$ is given by 
\begin{equation}\label{scal}
\scalar(g) = \unm\sum_k\tfrac{b_kd_k}{x_k}-\tfrac{1}{4}\sum_{i,j,k} \tfrac{x_k}{x_i x_j}[ijk].  
\end{equation}
We note that this formula is assuming nothing about $\Mm_{\mu_\pg}$, and if we set $Q:=-\kil_\ggo$ in the case when $G$ is semisimple, then $b_k=1$ for all $k$ and the formula is valid for any $\Ad(K)$-invariant decomposition as in \eqref{dec2}.  

A useful homothety invariant of $g$ is the {\it volume normalized scalar curvature} defined by 
\begin{equation}\label{scalN}
\scalar_N(g) := (\dete_Q{g})^{\frac{1}{n}}\scalar(g) = (x_1^{d_1}\dots x_r^{d_r})^{\frac{1}{n}}\scalar(g).  
\end{equation}
Note that $g\in\mca^G$ is Einstein if and only if $g$ is a critical point of $\scalar_N:\mca^G\rightarrow\RR$.  

If $G/K$ is not multiplicity-free, then it is sometimes useful to use instead formula \eqref{scal2} for the scalar curvature to find Einstein metrics (see \cite{Krr}).

\section{A formula for $\lic_\pg$ in terms of structural constants}\label{form-sec}

Let $g$ be a $G$-invariant metric on a connected homogeneous space $M=G/K$ with $G$ compact and $\ggo=\kg\oplus\pg$ a reductive decomposition as in \S\ref{str-sec}.  Recall that $\lic_\pg\sym(\pg)^H\subset\sym(\pg)^H$ for any intermediate subgroup $K\subset H\subset N_G(K)$, so 
$$
\Spec(\lic_\pg|_{\sym(\pg)^H})\subset\Spec(\lic_\pg|_{\sym(\pg)^K}).  
$$  
In the case when $\pg$ is $\Ad(H)$-multiplicity-free, we consider the orthonormal basis of $\sym(\pg)^H$ given by,
\begin{equation}\label{base}
\left\{ \tfrac{1}{\sqrt{d_1}}I_1,\dots, \tfrac{1}{\sqrt{d_r}}I_r\right\},
\end{equation}
where $I_k$ is the block map $\left[0,\dots,0,I_{\pg_k},0,\dots,0\right]$ of $\pg$.

\begin{theorem}\label{formLp}
Assume that the subspaces $\pg_k$'s are all $\Ad(H)$-invariant, $\Ad(H)$-irreducible and pairwise inequivalent for some intermediate subgroup $K\subset H\subset N_G(K)$.  Then, for any metric $g=(x_1,\dots,x_r)$, the entries of the matrix of the Lichnerowicz Laplacian $\lic_\pg$ with respect to the orthonormal basis \eqref{base} of $\sym(\pg)^H$ are given by
$$
[\lic_\pg]_{kk} = \tfrac{1}{d_k}\sum_{i,j\ne k} \tfrac{x_k}{x_ix_j}[ijk]  + \tfrac{1}{d_k}\sum_{i\ne k}  \tfrac{x_i}{x_k^2}[ikk], \qquad\forall k=1,\dots,r,
$$
and
$$
[\lic_\pg]_{km} = \tfrac{1}{\sqrt{d_k}\sqrt{d_m}} \sum_{i} \tfrac{x_i^2-x_k^2-x_m^2}{x_ix_kx_m} [ikm], \qquad\forall k\ne m.
$$
\end{theorem}

\begin{remark}\label{rem-Lp2}
In the case when $H=K$, i.e., $G/K$ is multiplicity-free, this is indeed the matrix of $\lic_\pg:\sym(\pg)^K\rightarrow\sym(\pg)^K$.   
\end{remark}

\begin{remark}\label{rem-Lp}
The spectrum of the symmetric $r\times r$ matrix given in the theorem, restricted to the hyperplane $\{(a_1,\dots,a_r):\sum \sqrt{d_i}a_i=0\}$ (which corresponds to $\sym_0(\pg)^H$), is therefore contained in $[\lambda_\pg, \lambda_\pg^{max}]$ (see \S\ref{Lp-sec}), providing estimates which are very useful in the study of the $G$-stability of an Einstein metric $g$.     
\end{remark}

\begin{remark}\label{Lp-gen}
If the decomposition $\pg=\pg_1\oplus\dots\oplus\pg_r$ is only assumed to be $\Ad(K)$-invariant and we call $L$ the symmetric $r\times r$ matrix defined as in the above theorem using the corresponding structural constants, then $L$ is a principal submatrix of the matrix of $\lic_\pg|_{\sym(\pg)^K}$ and so the spectrum of $L$ restricted to the hyperplane $\sum \sqrt{d_i}a_i=0$ (intersected with $\tca\tca_g^G$ if necessary) lies in $[\lambda_\pg, \lambda_\pg^{max}]$.  In particular, the $G$-instability of an Einstein metric $g$ (say $\ricci(g)=\rho g$) follows as soon as some eigenvalue of $L$ is less than $2\rho$.    
\end{remark}

\begin{proof}
Recall from \S\ref{ric-sec} that $g$ defines inner products on $\sym(\pg)$ and $\Lambda^2\pg^*\otimes\pg$ both denoted by $\ip$.  According to the definition of the operator $\lic_\pg$ given in \eqref{Lp},
\begin{equation}\label{Lpkm}
\la\lic_\pg I_k,I_m\ra = \unm\la\theta(I_m)\theta(I_k)\mu_\pg,\mu_\pg\ra + 2m_k\la I_k,I_m\ra, \qquad\forall k,m.
\end{equation}
If we write $\mu_\pg:\pg\times\pg\rightarrow\pg$ (see \eqref{mup}) as
$$
\mu_\pg = \sum_{i,j,l} \mu_{ij}^l, \qquad \mu_{ij}^l:\pg_i\times\pg_j\longrightarrow\pg_l,
$$
then following properties can be easily checked:
\begin{enumerate}[{\small $\bullet$}]
\item $\mu_{ij}^l(X,Y)=-\mu_{ji}^l(Y,X)$ for all $X\in\pg_i$, $Y\in\pg_j$,

\item $\mu_{ij}^l\perp\mu_{i'j'}^{l'}$ for any $(i,j,l)\ne (i',j',l')$,

\item $|\mu_{ij}^l|^2=\frac{x_l}{x_ix_j}[ijl]$; in particular, $\sum\limits_{i,j,k} \frac{x_k}{x_ix_j}[ijk] = |\mu_\pg|^2$.
\end{enumerate}

Using \eqref{tita} we obtain that
$$
\theta(I_k)\mu_{ij}^l = \left\{\begin{array}{cl}
0, & k\ne i,j,l,  \\
-\mu_{kj}^l, & k=i\ne j,l,\\
-\mu_{ik}^l, & k=j\ne i,l,\\
-2\mu_{kk}^l, & k=i=j\ne l,\\
\mu_{ij}^k, & k=l\ne i,j,\\
0, & k=l=i\ne j,\\
0, & k=l=j\ne i,\\
-\mu_{kk}^k, & k=l=i=j,
\end{array}\right.,
$$
and
$$
\theta(I_m)\theta(I_k)\mu_{ij}^l = \left\{\begin{array}{cl}
\mu_{kj}^l, & m=k=i\ne j,l,\\
\mu_{ik}^l, & m=k=j\ne i,l,\\
4\mu_{kk}^l, & m=k=i=j\ne l,\\
\mu_{ij}^k, & m=k=l\ne i,j,\\
\mu_{kk}^k, & m=k=l=i=j,\\
\hline
\mu_{km}^l, & k=i\ne m=j\ne l\ne k,\\
-\mu_{kj}^m, & k=i\ne m=l\ne j\ne k,\\
0, & k=i\ne m=j=l,\\
\hline
\mu_{mk}^l, & k=j\ne m=i\ne l\ne k,\\
-\mu_{ik}^m, & k=j\ne m=l\ne i\ne k,\\
0, & k=j\ne m=i=l,\\
\hline
-2\mu_{kk}^m, & k=i=j\ne m=l,\\
\hline
-\mu_{mj}^k, & k=l\ne m=i\ne j\ne k,\\
-\mu_{im}^k, & k=l\ne m=j\ne i\ne k,\\
-2\mu_{mm}^k, & k=l\ne m=i=j,
\end{array}\right..
$$
This implies that if we set $[ijk]_g:=|\mu_{ij}^k|^2$, then
\begin{align*}
\la\theta(I_k)\theta(I_k)\mu_\pg,\mu_\pg\ra =& \sum_{j,l\ne k} [kjl]_g +\sum_{i,l\ne k} [kil]_g +4\sum_{l\ne k} [kkl]_g +\sum_{i,j\ne k} [ijk]_g +[kkk]_g, \\
=& 2\sum_{i,j\ne k} [kij]_g +4\sum_{i\ne k} [kki]_g +\sum_{i,j\ne k} [ijk]_g +[kkk]_g,
\end{align*}
and for $k\ne m$,
\begin{align*}
\la\theta(I_m)\theta(I_k)\mu_\pg,\mu_\pg\ra =& \sum_{l\ne k,m} [kml]_g -\sum_{j\ne k,m}[kjm]_g +\sum_{l\ne k,m} [kml]_g -\sum_{i\ne k,m}[kim]_g \\
& -2[kkm]_g -\sum_{j\ne k,m} [mjk]_g -\sum_{i\ne k,m} [mik]_g -2[mmk]_g \\
=& 2\sum_{i\ne k,m} [kmi]_g  -2\sum_{i\ne k,m}[kim]_g \\
& -2[kkm]_g  -2\sum_{i\ne k,m} [mik]_g -2[mmk]_g.
\end{align*}
It now follows from \eqref{Lpkm} and \eqref{mkstr} that in terms of $\left\{ \frac{1}{\sqrt{d_1}}I_1,\dots, \frac{1}{\sqrt{d_r}}I_r\right\}$,
\begin{align*}
[\lic_\pg]_{kk} =& \tfrac{1}{d_k}\la\lic_\pg I_k,I_k\ra = \tfrac{1}{2d_k}\la\theta(I_k)\theta(I_k)\mu_\pg,\mu_\pg\ra + 2m_k \\
=& \tfrac{1}{2d_k}\left(2\sum_{i,j\ne k} [kij]_g +4\sum_{i\ne k} [kki]_g +\sum_{i,j\ne k} [ijk]_g +[kkk]_g -2\sum_{i,j} [kij]_g +\sum_{i,j} [ijk]_g\right), \\
=& \tfrac{1}{2d_k}\left(2\sum_{i,j\ne k} [kij]_g +4\sum_{i\ne k} [kki]_g +2\sum_{i,j\ne k} [ijk]_g +2[kkk]_g -2\sum_{i,j} [kij]_g +2\sum_{i\ne k} [ikk]_g\right), \\
=& \tfrac{1}{2d_k}\left(-2\sum_{i\ne k} [kki]_g -2\sum_{i\ne k} [kik]_g +4\sum_{i\ne k} [kki]_g +2\sum_{i,j\ne k} [ijk]_g  +2\sum_{i\ne k} [ikk]_g\right), \\
=& \tfrac{1}{d_k}\left(\sum_{i\ne k} [kki]_g +\sum_{i,j\ne k} [ijk]_g  \right),
\end{align*}
and for $k\ne m$,
\begin{align*}
[\lic_\pg]_{km} =& \tfrac{1}{\sqrt{d_k}\sqrt{d_m}}\la\lic_\pg I_m,I_k\ra = \tfrac{1}{2\sqrt{d_k}\sqrt{d_m}}\la\theta(I_m)\theta(I_k)\mu_\pg,\mu_\pg\ra  \\
=& \tfrac{1}{\sqrt{d_k}\sqrt{d_m}}\left( \sum_{i\ne k,m} [kmi]_g  -\sum_{i\ne k,m}[kim]_g
 -[kkm]_g  -\sum_{i\ne k,m} [mik]_g -[mmk]_g  \right).
\end{align*}
Finally, we use that $[ijk]_g=\frac{x_k}{x_ix_j}[ijk]$ to derive the formulas as stated in the theorem, concluding the proof.
\end{proof}

\begin{remark}
If $g$ is naturally reductive with respect to $G$ and $\pg$, then the numbers $[ijk]_g$ are invariant under any permutation of $ijk$ and so one obtains from the above proof that
$$
[\lic_\pg]_{kk} =
\tfrac{1}{d_k}\sum_{j\ne k;i} [ijk]_g, \quad\forall k,
\qquad
[\lic_\pg]_{km} =
-\tfrac{1}{\sqrt{d_k}\sqrt{d_m}}\sum_{i} [ikm]_g, \quad\forall k\ne m.
$$
This provides an alternative proof of \cite[Theorem 5.3]{stab-tres}.
\end{remark}

\begin{example}
Consider the full flag manifold $\SU(n)/T$, where $T$ is the diagonal maximal torus.  The standard reductive decomposition is given by 
$$
\sug(n)=\tg\oplus\pg_{12}\oplus\pg_{13}\oplus\dots\oplus\pg_{(n-1)n},
$$
where each $\pg_{ij}:=\{ zE_{ij}-\overline{z}E_{ji}:z\in\CC\}$ is $2$-dimensional and they are all $\Ad(T)$-irreducible and pairwise inequivalent.  Note that $\dim{\mca^G}=\frac{n(n-1)}{2}$.  Thus $[\pg_{ij},\pg_{kl}]_\pg=0$ if $\{ i,j\}$ and $\{ k,l\}$ are either equal or disjoint, and $[\pg_{ij},\pg_{ik}]_\pg$ is nonzero and it is contained in $\pg_{jk}$ for all $j\ne k$.  Moreover, since all the nonzero structural constants are equal to $\frac{1}{n}$, it follows from \eqref{rick} that the standard or Killing metric $g_{\kil}=(1,\dots,1)$ (set $Q=-\kil_{\sug(n)}$) is Einstein with $2\rho=1-\frac{1}{2n(n-2)}$.  According to Theorem \ref{formLp}, for $H=T$, one obtains that 
$
[\lic_\pg] = \frac{1}{2n}\Big(2(n-2)I-\Adj(X)\Big), 
$  
where $X=J(n,2,1)$ is the Johnson graph with parameters $(n,2,1)$ (see \cite[\S1.6]{GdsRyl}) and $\Adj(X)$ denotes its adjacency matrix.  The spectrum of $\Adj(X)$ is therefore given by 
$\{ 2(n-2), n-4, -2\}$ with multiplicities $1, n-1, \frac{n(n-3)}{2}$, respectively (see \cite[\S10.1, \S10.2]{GdsRyl}).  This implies that 
$$
\Spec(\lic_\pg)=\{0,\lambda_\pg,\lambda_\pg^{max}\}, \qquad
\lambda_\pg = \frac{1}{2}, \qquad \lambda_\pg^{max} = \frac{n-1}{n}, \qquad \forall n\geq 4,
$$  
with multiplicities $1$, $n-1$ and $\frac{n(n-3)}{2}$, respectively.  For $n=3$, $X$ is the complete graph on $3$ vertices and so the spectrum of $\Adj(X)$ equals $\{ 2, -1\}$, with respective multiplicities $1$ and $2$.  Thus $\lambda_\pg=\lambda_\pg^{max}=  \frac{1}{2}$ and has multiplicity $2$.  A straightforward comparison between the values of $2\rho$, $\lambda_\pg$ and $\lambda_\pg^{max}$ gives the following picture:  
\begin{enumerate}[{\small $\bullet$}] 
\item The standard metric $g_{\kil}$ on each $\SU(n)/T$ is always $G$-unstable with coindex $n-1$.  

\item $g_{\kil}$ is a local minimum for $n=3$, it is $G$-degenerate for $n=4$ and it is a saddle point for any $n\geq 5$.  

\item $g_{\kil}$ is always $G$-non-degenerate, and in particular $G$-rigid, except on $\SU(4)/T$.  We do not know whether $g_{\kil}$ is a local minimum for $\SU(4)/T$ or not.
\end{enumerate}
This example is a particular case of \cite[\S6]{stab-tres}, where the $G$-stability of the standard metric is studied using the spectrum of $\Adj(X)$ on $\SU(nk)/\Se(U(k)^n)$, $\Spe(nk)/\Spe(k)^n$ and $\SO(nk)/\Se(\Or(k)^n)$.  
\end{example}

\begin{example}\label{Stiefel}
The isotropy representation of the Stiefel manifold $\SO(k+2)/\SO(k)$, $k\geq 3$, decomposes as $\pg=\pg_1\oplus\pg_2\oplus\pg_3$, where $d_1=d_2=k$, $d_3=1$ and $\pg_1,\pg_2$ are equivalent $\Ad(K)$-representations of real type.  This implies that $\dim{\mca^G}=4$.  It is proved in \cite[\S 4]{Krr} that this space admits a unique Einstein metric given by $g=(1,1,\frac{2k}{k+1})$, with $2 \rho=\frac{k}{k+1}$, and that for the intermediate subgroup $H:=\SO(k)\times\SO(2)$ of the normalizer of $\SO(k)$ one has that $\mca^G=\Ad(H)\cdot\mca^{G,H}$.  Thus we only need to compute the spectrum of $\lic_\pg(g)$ restricted to $\sym(\pg)^H$ (see Remark \ref{H}), which is precisely the $3$-dimensional space of diagonal $\Ad(H)$-invariant tensors since $\pg_1,\pg_2$ are now inequivalent $\Ad(H)$-representations.  It follows from Theorem \ref{formLp} that
$$
L_\pg|_{\sym(\pg)^H} = \left[\begin{matrix}
\frac{k+1}{2k^2} &  \frac{k^2-2k-1}{2k^2(k+1)}  & -\frac{\sqrt{k}}{k+1} \\
\frac{k^2-2k-1}{2k^2(k+1)} &   \frac{k+1}{2k^2} & -\frac{\sqrt{k}}{k+1}\\
-\frac{\sqrt{k}}{k+1} & -\frac{\sqrt{k}}{k+1} &  \frac{2k}{k+1}
\end{matrix}\right], 
$$
and it is therefore straightforward to check that $\lambda_\pg=\tfrac{2k+1}{k^2(k+1)}$ and $\lambda_\pg^{max}=\tfrac{2k+1}{k+1}$.  Since $\lambda_\pg<2\rho<\lambda_\pg^{max}$, we obtain that $g$ is always a $G$-non-degenerate and $G$-unstable saddle point.    
\end{example}

\begin{example}\label{Stiefel3} 
In the case $S^2\times S^3=\SO(4)/\SO(2)$, one also has as above that $\pg=\pg_1\oplus\pg_2\oplus\pg_3$,  $d_1=d_2=2$, $d_3=1$ and $\pg_1,\pg_2$ equivalent, but here $\dim{\mca^G}=5$ since $\pg_1$ and $\pg_2$ are of complex type.  Besides the Einstein metric $g=(1,1,\frac{4}{3})$, there is another Einstein metric $g_\times\in\mca^G$, which was discovered in \cite{AlkFrrDtt} and is isometric to the product metric.  The metric $g_\times$ actually belongs to a one-parameter family of Einstein metrics in $\mca^G$ given by the orbit $\Aut(G/K)\cdot g_\times$ (note that $\Aut(G/K)=\Or(2)$).  The argument used in Example \ref{Stiefel} to compute the spectrum of $\lic_\pg(g)$ is still valid, giving that $\lambda_\pg(g)=\tfrac{5}{12}$ and $\lambda_\pg^{max}(g)=\tfrac{5}{3}$.  On the other hand, we have that $\dim{T_{g_\times}\Aut(G/K)\cdot g_\times}=1$ (see \eqref{tang}) and it is easy to see by using that $g_\times$ is the product metric that $g_\times$ is a $G$-non-degenerate and $G$-unstable saddle point of coindex $1$.    
\end{example}

Beyond left-invariant metrics on Lie groups, the above is the only example we know of a continuous family of Einstein metrics in $\mca_1^G$, i.e., with $\dim{T_{g}\Aut(G/K)\cdot g}>0$.

\subsection{Case $r=2$}\label{r2-sec}
For homogeneous spaces with only two isotropy summands the above formulas considerably simplify.  We refer to \cite{DckKrr, He} for a complete classification of these spaces  ($43$ infinite families and $78$ isolated examples) together with their Einstein metrics (only $19$ isolated examples do not admit any).  Note that $\dim{\mca_1^G}=1$ in this case (with the only exception of $\SO(8)/G_2$), so $\scalar|_{\mca_1^G}$ is a one-variable function.

It follows from \eqref{rick} that the Ricci eigenvalues are
$$
\begin{array}{l}
\rho_1 = \left(\tfrac{b_1}{2} -\tfrac{[111]}{4d_1} -\tfrac{[122]}{2d_1}\right)\tfrac{1}{x_1} -\tfrac{[112]}{2d_1}\tfrac{x_2}{x_1^2} + \tfrac{[122]}{4d_1}\tfrac{x_1}{x_2^2}, 
 \\ 
\rho_2 = \left(\tfrac{b_2}{2} -\tfrac{[222]}{4d_2} -\tfrac{[112]}{2d_2}\right)\tfrac{1}{x_2} -\tfrac{[122]}{2d_2}\tfrac{x_1}{x_2^2} + \tfrac{[112]}{4d_2}\tfrac{x_2}{x_1^2}, 
\end{array}
$$
and from Theorem \ref{formLp} that
$$
[\lic_\pg]_{11} = \tfrac{[122]x_1}{d_1x_2^2} + \tfrac{[112]x_2}{d_1x_1^2}, \qquad 
[\lic_\pg]_{22} = \tfrac{[112]x_2}{d_2x_1^2} + \tfrac{[122]x_1}{d_2x_2^2}.   
$$
Since $\Spec(\lic_\pg)=\{ 0,\lambda_\pg\}$, we obtain that
$$
\lambda_\pg = [\lic_\pg]_{11} + [\lic_\pg]_{22} = \tfrac{d_1+d_2}{d_1d_2} \left(\tfrac{[112]x_2}{x_1^2} + \tfrac{[122]x_1}{x_2^2}\right).   
$$
This can be used to find the $G$-stability and critical point types of any Einstein metric $(x_1,x_2)$ on a homogeneous space $M=G/K$ of this kind, after computing or finding in the literature the structural constants of $G/K$ (see \cite{stab} for the case of standard Einstein metrics).   

\begin{example}\label{r2-flag}
If $M=G/K$ is a flag manifold (i.e., $G$ simple and $K$ the centralizer of a torus in $G$, see \S\ref{flag-sec}) with two isotropy summands (see e.g.\ \cite[Table 1]{AnsChr} for a complete list, consisting of $3$ classical families and $10$ exceptional examples), then it can be assumed that the only non-zero structural constant is $[112]=\frac{d_1d_2}{d_1+4d_2}$ (set $Q:=-\kil_{\ggo}$) and so, 
$$
\rho_1 = \tfrac{1}{2x_1} -\tfrac{[112]x_2}{2d_1x_1^2}, \qquad 
\rho_2 = \tfrac{1}{2x_2} - \tfrac{[112](2x_1^2-x_2^2)}{4d_2x_1^2x_2}, \qquad
\lambda_\pg = \tfrac{(d_1+d_2)x_2}{(d_1+4d_2)x_1^2}.  
$$
It is easy to see that there are exactly two solutions to $\rho_1=\rho_2$ up to scaling, the K\"ahler-Einstein metric $g_0=(1,2)$ and a non-K\"ahler Einstein metric given by $g_1=(1,\frac{4d_2}{d_1+2d_2})$.  On the other hand, a straightforward inspection gives that $g_0$ is always $G$-stable (i.e., $2\rho<\lambda_\pg$) and $g_1$ is always $G$-unstable (i.e., $\lambda_\pg<2\rho$).  In particular, $g_0$ is a local maximum of $\scalar|_{\mca_1^G}$ and $g_1$ a local minimum for any of these spaces (cf.\ \cite[Theorem 1]{AnsChr2}).      
\end{example}

\begin{table}{\small 	
$$
\begin{array}{c|c|c|c|c|c}
W_i &\ggo/\kg & a_1=\frac{[123]}{d_1} & a_2=\frac{[123]}{d_2} & a_3=\frac{[123]}{d_3} & [123] \\
[2mm]  \hline \hline
\rule{0pt}{14pt}
W_1&\frac{\sog(k+l+m)}{\sog(k)\oplus\sog(l)\oplus\sog(m)} &\frac{k}{2(k+l+m-2)} & \frac{l}{2(k+l+m-2)} & \frac{m}{2(k+l+m-2)} & \frac{klm}{2(k+l+m-2)}
\\[2mm] \hline\rule{0pt}{14pt}
W_2&\frac{\sug(k+l+m)}{\sg(\ug(k)\oplus\ug(l)\oplus\ug(m))} & \frac{k}{2(k+l+m)} &\frac{l}{2(k+l+m)} & \frac{m}{2(k+l+m)} &  \frac{klm}{k+l+m}
\\[2mm] \hline\rule{0pt}{14pt}
W_3&\frac{\spg(k+l+m)}{\spg(k)\oplus\spg(l)\oplus\spg(m)} &\frac{k}{2(k+l+m+1)} &\frac{l}{2(k+l+m+1)} & \frac{m}{2(k+l+m+1)} & \frac{2klm}{k+l+m+1}
\\[2mm] \hline\rule{0pt}{14pt}
W_4&\frac{\sug(2l)}{\ug(l)}, l\geq 2 & \frac{l^2-1}{4} & \frac{l-1}{4l} & \frac{1}{4} &  \frac{l+1}{4l}
\\[2mm] \hline \rule{0pt}{14pt}
W_5&\frac{\sog(2l)}{\ug(1)\oplus\ug(l-1)},\,l\geq 4  & \frac{l-2}{4(l-1)} & \frac{l-2}{4(l-1)} & \frac{1}{2(l-1)} & \frac{l-2}{2}
\\[2mm] \hline \rule{0pt}{14pt}
W_6&\frac{\eg_6}{\sug(4)\oplus 2\spg(1)\oplus\RR} & \frac{1}{4} &\frac{1}{4} & \frac{1}{6} & 4
\\[2mm] \hline \rule{0pt}{14pt}
W_7&\frac{\eg_6}{\sog(8) \oplus \RR^2}  & \frac{1}{6} & \frac{1}{6} & \frac{1}{6} & \frac{8}{3}
\\[2mm] \hline \rule{0pt}{14pt}
W_8&\frac{\eg_6}{\spg(3)\oplus\spg(1)} & \frac{1}{4} & \frac{1}{8} & \frac{7}{24} & \frac{7}{2}
\\[2mm] \hline \rule{0pt}{14pt}
W_9&\frac{\eg_7}{\sog(8)\oplus 3\spg(1)} &\frac{2}{9} & \frac{2}{9} &\frac{2}{9} & \frac{64}{9}
\\[2mm] \hline \rule{0pt}{14pt}
W_{10}&\frac{\eg_7}{\sug(6)\oplus \spg(1) \oplus \RR} & \frac{2}{9} & \frac{1}{6} & \frac{5}{18} & \frac{20}{3}
\\[2mm] \hline \rule{0pt}{14pt}
W_{11}&\frac{\eg_7}{\sog(8)}  & \frac{5}{18} &\frac{5}{18} & \frac{5}{18} & \frac{175}{18}
\\[2mm] \hline \rule{0pt}{14pt}
W_{12}&\frac{\eg_8}{\sog(12)\oplus 2\spg(1)} & \frac{1}{5} &\frac{1}{5} & \frac{4}{15} & \frac{64}{5}
\\[2mm]  \hline \rule{0pt}{14pt}
W_{13}&\frac{\eg_8}{\sog(8)\oplus \sog(8)}  & \frac{4}{15} & \frac{4}{15} & \frac{4}{15} & \frac{256}{15}
\\[2mm]  \hline \rule{0pt}{14pt}
W_{14}&\frac{\fg_4}{\sog(5)\oplus 2\spg(1)}  & \frac{5}{18} & \frac{5}{18} & \frac{1}{9} & \frac{20}{9}
\\[2mm] \hline \rule{0pt}{14pt}
W_{15}&\frac{\fg_4}{\sog(8)}  & \frac{1}{9} & \frac{1}{9} & \frac{1}{9} & \frac{8}{9}
\\[2mm] \hline \hline
\end{array}
$$}
\caption{Generalized Wallach spaces with $G$ simple (see \cite{Nkn, ChnKngLng}).  In case $W_1$, the triple $(k,2,2)$ (and its permutations) must be avoided, and $(k,1,1)$ (i.e., the Stiefel manifold $\SO(k+2)/\SO(k)$) are the only instances in which $G/K$ is not multiplicity-free.  The spaces $W_2$, $W_5$ and $W_7$ are also flag manifolds (see \S\ref{flag-sec}).}
\label{W1}
\end{table}

\begin{table}{\small 	
$$
\begin{array}{c|c|c|c|c|c|c|c|c}
W_i &\ggo/\kg & b &\frac{1-2b}{2b} & g &\lambda_\pg&\lambda_\pg^{max}&2\rho & \text{Type}\\
[2mm]  \hline \hline
\rule{0pt}{14pt}
W_{1}&\sog(3) &\unm& & g_{\kil} & \frac{3}{2} & \frac{3}{2} & \frac{1}{2}& {\bf G\text{{\bf -stab}.}}
\\[2mm] \hline\rule{0pt}{14pt}
\multirow{2}{*}{$W_{1}$} & \frac{\sog(3k)}{\sog(k)\oplus\sog(k)\oplus\sog(k)} &\multirow{2}{*}{$\frac{k}{2(3k-2)}$} & \multirow{2}{*}{$\frac{2(k-1)}{k}$} &
g_{\kil} &\frac{3k}{2(3k-2)}&\frac{3k}{2(3k-2)}&\frac{5k-4}{2(3k-2)} & G\text{-unst.} \\ [2mm]
& k \ge 3 & & & g_i & \frac{2k-1}{(k-1)(3k-2)} & \frac{3(k-1)}{3k-2} & \frac{2k-1}{3k-2} & G\text{-unst.}
\\[2mm] \hline\rule{0pt}{14pt}
\multirow{2}{*}{$W_2$}&\multirow{2}{*}{$\frac{\sug(3k)}{\sg(\ug(k)\oplus\ug(k)\oplus\ug(k))}$} & \multirow{2}{*}{$\frac{1}{6}$} & \multirow{2}{*}{$2$} & g_{\kil} & \unm & \unm &  \frac{5}{6} & G\text{-unst.} \\ [2mm]
 & & && g_i & 0 & 1 & \frac{2}{3}& G\text{-unst.}
\\[2mm] \hline\rule{0pt}{14pt}
\multirow{2}{*}{$W_3$}&\multirow{2}{*}{$\frac{\spg(3k)}{\spg(k)\oplus\spg(k)\oplus\spg(k)}$} & \multirow{2}{*}{$\frac{k}{2(3k+1)}$} & \multirow{2}{*}{$\frac{2k+1}{k}$} &
 g_{\kil}, & \frac{3k}{2(3k+1)} &\frac{3k}{2(3k+1)} & \frac{5k+2}{2(3k+1)} & G\text{-unst.} \\ [2mm]
& & && g_i  & \frac{-(4k+1)}{2(2k+1)(3k+1)} & \frac{3(2k+1)}{2(3k+1)} & \frac{4k+1}{2(3k+1)}&G\text{-unst.}
\\[2mm] \hline\rule{0pt}{14pt}
\multirow{2}{*}{$W_5$}&\multirow{2}{*}{$\frac{\sog(8)}{\ug(1)\oplus\ug(3)}$}  & \multirow{2}{*}{$\frac{1}{6}$} & \multirow{2}{*}{$2$} &
  g_{\kil} & \unm & \unm &  \frac{5}{6} & G\text{-unst.} \\ [2mm]
 & & && g_i & 0 & 1 & \frac{2}{3} & G\text{-unst.}
 \\[2mm] \hline \rule{0pt}{14pt}
\multirow{2}{*}{$W_7$}&\multirow{2}{*}{$\frac{\eg_6}{\sog(8) \oplus \RR^2}$}  & \multirow{2}{*}{$\frac{1}{6}$} & \multirow{2}{*}{$2$} &
  g_{\kil} & \unm & \unm &  \frac{5}{6} & G\text{-unst.} \\ [2mm]
 & && & g_i & 0 & 1 & \frac{2}{3} & G\text{-unst.}
\\[2mm] \hline \rule{0pt}{14pt}
\multirow{2}{*}{$W_9$}&\multirow{2}{*}{$\frac{\eg_7}{\sog(8)\oplus 3\spg(1)}$} & \multirow{2}{*}{$\frac{2}{9}$} & \multirow{2}{*}{$\frac{5}{4}$} &
 g_{\kil} & \frac{2}{3} & \frac{2}{3} & \frac{7}{9} & G\text{-unst.} \\ [2mm]
 & && &  g_i & \frac{13}{30} & \frac{5}{6} & \frac{13}{18} & G\text{-unst.}
\\[2mm] \hline \rule{0pt}{14pt}
\multirow{2}{*}{$W_{11}$}&\multirow{2}{*}{$\frac{\eg_7}{\sog(8)}$} & \multirow{2}{*}{$\frac{5}{18}$} & \multirow{2}{*}{$\frac{4}{5}$} &
 g_{\kil} & \frac{15}{18} &\frac{15}{18} &  \frac{13}{18} & {\bf G\text{{\bf -stab}.}} \\ [2mm]
 & & &&  g_i &  \frac{2}{3} & \frac{7}{6} & \frac{7}{9} & G\text{-unst.}
\\[2mm] \hline \rule{0pt}{14pt}
\multirow{2}{*}{$W_{13}$}& \multirow{2}{*}{$\frac{\eg_8}{\sog(8)\oplus \sog(8)}$}  & \multirow{2}{*}{$\frac{4}{15}$} & \multirow{2}{*}{$\frac{7}{8}$} &
 g_{\kil} &\frac{12}{15} & \frac{12}{15} & \frac{11}{15} & {\bf G\text{{\bf -stab}.}} \\ [2mm]
& & &&  g_i  & \frac{7}{10}  & \frac{69}{70} & \frac{23}{30} & G\text{-unst.}
\\[2mm]  \hline \rule{0pt}{14pt}
\multirow{2}{*}{$W_{15}$}&\multirow{2}{*}{$\frac{\fg_4}{\sog(8)}$} & \multirow{2}{*}{$\frac{1}{9}$} & \multirow{2}{*}{$\frac{7}{2}$} &
 g_{\kil} & \frac{1}{3} & \frac{1}{3} & \frac{8}{9} & G\text{-unst.} \\ [2mm]
& & & &g_i  & -\frac{11}{42} & \frac{7}{6} & \frac{11}{18} & G\text{-unst.}
\\[2mm] \hline \hline
\end{array}
$$}
\caption{Einstein metrics $g_{\kil},g_1,g_2,g_3$ on generalized Wallach spaces with $a_1=a_2=a_3=b$ (see \cite{FrsLmsNkn}).  The metrics $g_i$'s all satisfy that $\lambda_\pg<2\rho<\lambda_\pg^{max}$ and so they are $G$-unstable and $G$-non-degenerate saddle points with coindex $1$.  All the standard metrics $g_{\kil}$'s which are $G$-unstable are {\bf local minima} of coindex $2$ since $\lambda_\pg=\lambda_\pg^{max}<2\rho$.  Note that $g_1,g_2,g_3$ are all K\"ahler metrics in the flag cases $W_2$, $W_5$ and $W_7$, which are precisely the cases when $\lambda_\pg=0$.}
\label{W2}
\end{table}

\begin{table}{\small 	
$$
\begin{array}{c|c|c||c|c|c}
W_i  & g & (\det_{g_{\kil}}g)^{\frac{1}{n}}\scalar(g) & W_i  & g & (\det_{g_{\kil}}g)^{\frac{1}{n}}\scalar(g)\\
[2mm]  \hline \hline
\rule{0pt}{14pt}
\multirow{2}{*}{$W_{1}$} & \multirow{2}{*}{$g_{\kil}$} & \multirow{2}{*}{$\frac{3}{4} \quad (k=3)$} & \multirow{2}{*}{$W_{7}$} & g_{\kil} & 20  \\ [2mm]
& & & & g_i & 16 \cdot 2^\frac{1}{3} \approx 20.1587
\\[2mm] \hline\rule{0pt}{14pt}
\multirow{2}{*}{$W_{1}$} & g_{\kil} & \frac{3k^2(5k-4)}{4(3k-2)} & \multirow{2}{*}{$W_{9}$} & g_{\kil} & \frac{112}{3} \approx 37.3333  \\ [2mm]
& g_i & \frac{3k^2(2k-1)}{3k-2}(\frac{2k-2}{k})^{\frac{1}{3}} & & g_i & \frac{52}{3} {10}^\frac{1}{3} \approx 37.3435
\\[2mm] \hline\rule{0pt}{14pt}
\multirow{2}{*}{$W_{2}$} & g_{\kil} &  \frac{5k^2}{2} & \multirow{2}{*}{$W_{11}$} & g_{\kil} & \frac{455}{12} \approx 37.9166  \\ [2mm]
& g_i &  2^{\frac{4}{3}} k^2 \approx 2.52 k^2 & & g_i & \frac{49}{6} {10}^\frac{2}{3} \approx 37.9063
\\[2mm] \hline\rule{0pt}{14pt}
\multirow{2}{*}{$W_{3}$} & g_{\kil} & \frac{n(5k+2)}{4(3k+1)} & \multirow{2}{*}{$W_{13}$} & g_{\kil} & \frac{352}{5} = 70.4  \\ [2mm]
& g_i & \frac{3k^2(4k+1)}{3k+1} (\frac{2*k+1}{k})^\frac{1}{3} & & g_i & \frac{184}{5} {7}^\frac{2}{3} \approx 70.3958
\\[2mm] \hline\rule{0pt}{14pt}
\multirow{2}{*}{$W_{5}$} & g_{\kil} & \frac{15}{2}=7.5 & \multirow{2}{*}{$W_{15}$} & g_{\kil} & \frac{32}{3} \approx 10.6666  \\ [2mm]
& g_i & 6 \cdot 2^\frac{1}{3} \approx 7.5595 & & g_i & \frac{11}{3} {28}^\frac{1}{3} \approx 11.1341
\\[2mm] \hline \hline
\end{array}
$$}
\caption{Volume normalized scalar curvature of Einstein metrics in Table \ref{W2}.  Note that $\scalar_N(g_{\kil}) < \scalar_N(g_i)$ always holds, except for $W_{11}$ and $W_{13}$, that is, precisely when the standard metric $g_{\kil}$ is a local maximum.}
\label{W2Sc}
\end{table}

\section{Generalized Wallach spaces}\label{W-sec}

A connected homogeneous space $M=G/K$ with $G$ compact semisimple is called a {\it generalized Wallach space} if the $\kil_\ggo$-orthogonal reductive decomposition $\ggo=\kg\oplus\pg$ satisfies that $\pg=\pg_1\oplus\pg_2\oplus\pg_3$ for some $\Ad(K)$-irreducible subspaces $\pg_i$'s such that $[\pg_i,\pg_i]\subset\kg$ for all $i=1,2,3$.  Equivalently, the only nonzero structural constant is $[123]$.  The generalized Wallach spaces with $G$ simple are listed in Table \ref{W1}, together with the numbers $[123]$ and $a_i:=[123]/d_i$, where $d_i:=\dim{\pg_i}$ (see \cite{Nkn, ChnKngLng}).  For a table with the dimensions $d_i$'s see \cite[Table 1]{ChnNkn}.  

In what follows, as an application of the formula given in Theorem \ref{formLp}, we compute the spectrum of the Lichnerowicz Laplacian $\lic_\pg$ for each $G$-invariant Einstein metric on most of these spaces (see Tables \ref{W2}, \ref{W4} and \ref{W5}).  In this way, we obtain the $G$-stability type of all these Einstein metrics, as well as what type of critical points of $\scalar|_{\mca_1^G}$ are.  It was proved in \cite[Theorem 7]{AbvArvNknSss} that, generically, there are three possible behaviors depending on the region the triple $(a_1,a_2,a_3)$ attached to $G/K$ belongs.  The regions are described in \cite[Table 2]{ChnNkn} (the case $5$ is incorrectly assigned to region $O_3$, it is actually in $O_1$).  It is easy to check that the spaces $W_1$ with $l=m$ are in $O_1$.  There is no contradiction between our results and the behaviors described in these papers.    

Note that $\dim{\mca_1^G}=2$ for any generalized Wallach space (except $\SO(k+2)/\SO(k)$, see Example \ref{Stiefel}) and so
$$
\Spec(\lic_\pg(g))=\{ 0,\lambda_\pg,\lambda_\pg^{max}\}
$$
for any Einstein metric $g=(x_1,x_2,x_3)\in\mca^G$.  We will always set $Q=-\kil_\ggo$ as a background metric.  Using that $[123]$ is the only nonzero structural constant, it easily follows from Theorem \ref{formLp} that 
\begin{equation}\label{Lpgen}
\lic_\pg(g) = 
\left[\begin{matrix}
 a_1\frac{2x_1}{x_2x_3} & \sqrt{a_1a_2}\tfrac{x_3^2-x_1^2-x_2^2}{x_1x_2x_3} & \sqrt{a_1a_3}\tfrac{x_2^2-x_1^2-x_3^2}{x_1x_2x_3}\\
  \sqrt{a_1a_2}\tfrac{x_3^2-x_1^2-x_2^2}{x_1x_2x_3} & a_2\frac{2x_2}{x_1x_3} &\sqrt{a_2a_3}\tfrac{x_1^2-x_2^2-x_3^2}{x_1x_2x_3}\\
  \sqrt{a_1a_3}\tfrac{x_2^2-x_1^2-x_3^2}{x_1x_2x_3} & \sqrt{a_2a_3}\tfrac{x_1^2-x_2^2-x_3^2}{x_1x_2x_3} & a_3\frac{2x_3}{x_1x_2} 
\end{matrix}\right].  
\end{equation}

\subsection{Case $d_1=d_2=d_3$}\label{igu-sec}
It is proved in \cite[Theorem 3]{FrsLmsNkn} that $G/K$ ($\ne\SO(3)$) admits in this case four $G$-invariant Einstein metrics up to scaling, the standard metric $g_{\kil}=(1,1,1)$ and
$$
g_1=(\tfrac{1-2b}{2b},1,1), \qquad g_2=(1,\tfrac{1-2b}{2b},1), \qquad g_3=(1,1,\tfrac{1-2b}{2b}),
$$
where $b:=a_1=a_2=a_3=[123]/d_i$.  Note that $b\leq\unm$ and equality holds if and only if $\ggo/\kg=\sog(3)$.  It follows from \eqref{rick} that for $g_{\kil}$, $2 \rho= 1-b$ and by \eqref{Lpgen},
$$
\lic_\pg(g_{\kil})= \left[\begin{matrix}
2b & - b  & - b \\
- b & 2b & -b\\
- b & - b & 2b
\end{matrix}\right], \qquad \mbox{so}\quad \lambda_\pg=\lambda_\pg^{max}=3b.
$$
On the other hand, for each $g_i$ we have that $2 \rho= \frac{1+2b}{2}$,
$$
L_\pg(g_1) = \left[\begin{matrix}
(1-2b)& - \unm (1-2b) & - \unm (1-2b) \\
- \unm (1-2b)& \tfrac{4b^2}{1-2b} & -\tfrac{4b^2+4b-1}{2(1-2b)} \\
- \unm (1-2b) & -\tfrac{4b^2+4b-1}{2(1-2b)}  &  \tfrac{4b^2}{1-2b}
\end{matrix}\right],
$$
and thus it is straightforward to check that
$$
\Spec(\lic_\pg(g_i))=\left\{ 0, \tfrac{12b^2+4b-1}{2(1-2b)}, \tfrac{3(1-2b)}{2}\right\}, \qquad\forall i=1,2,3.
$$
All this data together with the resulting types have been collected in Table \ref{W2}.  We furthermore give in Table \ref{W2Sc} the value of the volume normalized scalar curvature $\scalar_N$ (see \eqref{scalN}) for these metrics.  

The two $G$-stable standard metrics found on $E_7/\SO(8)$ and $E_8/(\Spin(8)\times\Spin(8))$ deserve special attention, note that they are in particular local maxima of $\scalar|_{\mca_1^G}$.  We do not know  whether these Einstein metrics are stable in the classical sense (see Remark \ref{new}).  We do not know whether they realize the Yamabe invariant either.

\subsection{Case $d_1=d_2\ne d_3$}\label{dist-sec}
We set $b:=a_1=a_2$ and $c:=a_3$.  According to \cite[Theorem 4]{FrsLmsNkn}, under this assumption, $G/K$ admits at most four $G$-invariant Einstein metrics given by:

\begin{enumerate}[{\small $\bullet$}]
\item $g_{q^\pm}=(1,1,q^\pm)$, where
$$
q^\pm= \tfrac{1 \pm \sqrt{1-4(b+c)(1-2c)}}{2(b+c)},
$$
\item $g_{p^\pm}=(p^\pm,1,2b(p^\pm + 1))$, where
$$
p^\pm= \tfrac{-1+2b-8b^2(b+c) \pm \sqrt{(1-2b+8b^2(b+c))^2-4(b+c)^2(1-4b^2)^2}}{2(b+c)(1-4b^2)}.
$$
Note that $p^+p^-=1$ and $1<p^+$.  The radicand is positive if and only if
\begin{equation}\label{T}
T:=1-2(2b+c)+16b^2(b+c)>0.
\end{equation}
\end{enumerate}
It is worth pointing out that the non-negativity of the radicands is a necessary condition for the existence of the corresponding metric, which may o may not hold.

\begin{table}{\small 	
$$
\begin{array}{c|c|c|c|c}
W_i & \ggo/\kg & b & c & g \\
[2mm]  \hline \hline
\rule{0pt}{14pt}
\multirow{2}{*}{$W_1$} & \multirow{2}{*}{$\frac{\sog(4)}{\sog(2)}$} & \multirow{2}{*}{$\frac{1}{4}$}  & \multirow{2}{*}{$\unm$}  & q^+ = \frac{4}{3}\\ [2mm]
& & & & g_{\times}  \; (\text{see Ex.\ \ref{Stiefel3}})
\\[2mm] \hline\rule{0pt}{14pt}
\multirow{2}{*}{$W_1$} & \frac{\sog(k+2)}{\sog(k)} & \multirow{2}{*}{$\frac{1}{2k}$}  & \multirow{2}{*}{$\unm$}  & \multirow{2}{*}{$q^+ = \frac{2k}{k+1} \; (\text{see Ex.\ \ref{Stiefel}})$}\\ [2mm]
& k\geq 3 & & & 
\\[2mm] \hline\rule{0pt}{14pt}
\multirow{2}{*}{$W_1$} & \frac{\sog(k+2m)}{\sog(k)\oplus2\sog(m)} & \multirow{2}{*}{$\frac{m}{2(k+2m-2)}$} & \multirow{2}{*}{$\frac{k}{2(k+2m-2)}$} & q^\pm=\frac{k+2m-2 \pm \sqrt{k^2-4m+4}}{k+m}\\ [2mm]
& m\geq 3 & & & p^\pm  \; (\text{see \eqref{qpFam1}})
\\[2mm] \hline\rule{0pt}{14pt}
\multirow{4}{*}{$W_2$}&\multirow{4}{*}{$\frac{\sug(k+2m)}{\sg(\ug(k)\oplus 2\ug(m))}$} & \multirow{4}{*}{$\frac{m}{2(k+2m)}$} &
\multirow{4}{*}{$\frac{k}{2(k+2m)}$} & q^+=2 \\ [2mm]
& & & & q^-=\frac{2m}{k+m}
\\ [2mm]
& & & & p^+=\frac{k+3m}{k+m}   \\ [2mm]
& & & & p^-=\frac{k+m}{k+3m}
\\[2mm] \hline\rule{0pt}{14pt}
\multirow{2}{*}{$W_3$}&\multirow{2}{*}{$\frac{\spg(k+2m)}{\spg(k)\oplus 2\spg(m)}$} & \multirow{2}{*}{$\frac{m}{2(k+2m+1)}$} &
\multirow{2}{*}{$\frac{k}{2(k+2m+1)}$}& q^\pm=\frac{k+2m+1\pm\sqrt{k^2+2m+1}}{k+m}
\\ [2mm]
& & & &  p^\pm \; (\text{see \eqref{pmmFam3}})
\\[2mm] \hline\rule{0pt}{14pt}
\multirow{4}{*}{$W_5$}&\multirow{4}{*}{$\frac{\sog(2l)}{\ug(1)\oplus\ug(l-1)}$}  & \multirow{4}{*}{$\frac{l-2}{4(l-1)}$} &
 \multirow{4}{*}{$\frac{1}{2(l-1)}$} & q^+=2 \\ [2mm]
 &  & & &  q^-=\frac{2(l-2)}{l}
 \\ [2mm]
 & l\geq 5 & & &  p^+=\frac{3l-4}{l}  \\ [2mm]
 & & & &  p^-=\frac{l}{3l-4}
 \\[2mm] \hline \rule{0pt}{14pt}
\multirow{2}{*}{$W_6$}&\multirow{2}{*}{$\frac{\eg_6}{\sug(4) \oplus 2\spg(1) \oplus \RR}$}  & \multirow{2}{*}{$\frac{1}{4}$} &
  \multirow{2}{*}{$\frac{1}{6}$} &  p^+= \frac{5}{3} \\ [2mm]
 & & & & p^-=\frac{3}{5}
\\[2mm] \hline \rule{0pt}{14pt}
\multirow{2}{*}{$W_{12}$}&\multirow{2}{*}{$\frac{\eg_8}{\sog(12)\oplus 2\spg(1)}$} & \multirow{2}{*}{$\frac{1}{5}$} & \multirow{2}{*}{$\frac{4}{15}$} &  q^+= \frac{15+\sqrt{29}}{14}  \\ [2mm]
 & & & &  q^-=\frac{15-\sqrt{29}}{14}
\\[2mm] \hline \rule{0pt}{14pt}
\multirow{2}{*}{$W_{14}$}&\multirow{2}{*}{$\frac{\fg_4}{\sog(5)\oplus2\spg(1)}$} & \multirow{2}{*}{$\frac{5}{18}$} &
  \multirow{2}{*}{$\frac{1}{9}$}& p^+=\frac{499+9\sqrt{1177}}{392}  \\ [2mm]
 & & & & p^-=\frac{499-9\sqrt{1177}}{392}
\\[2mm] \hline \hline
\end{array}
$$}
\caption{Einstein metrics on generalized Wallach spaces with $a_1=a_2=b\ne a_3=c$, given by $g_{q^\pm}=(1,1,q^{\pm})$ and $g_{p^\pm}=(p^\pm,1,2b(p^\pm +1))$ (see \cite{FrsLmsNkn}).  Assume that $k\ne m$ everywhere.}
\label{W3}
\end{table}

\begin{table}{\small 	
$$
\begin{array}{c|c|c}
W_i & g & (\det_{g_{\kil}}g)^{\frac{1}{n}}\scalar(g)
\\
[2mm]  \hline \hline
\rule{0pt}{14pt}
%
%
\multirow{2}{*}{$W_1$} & q^\pm & \frac{m(2k+m)(4m+2k-4-m q)}{4(2m+k-2)} q^\frac{m}{2k+m}\\ [2mm]
m\geq 3 & p^\pm & \frac{(p+1)m(2k+ m)(3m^2+4km-8m+k^2-4k+4)}{4p(2m+k-2)^2} p^{\frac{m}{2k+m}} \left(\frac{m(p+1)}{2m+k-2}\right)^{\frac{k}{2km+m^2}}
\\[2mm] \hline\rule{0pt}{14pt}
\multirow{4}{*}{$W_2$} & q^+ &  \frac{m(2k+m)(m+k)}{(2m+k)} 2^{\frac{m}{2k+m}} \\ [2mm]
 & q^- & \frac{m(2k+m)(m^2+3km+k^2)}{(2m+k)(m+k)} \left(\frac{2m}{m+k}\right)^\frac{m}{2k+m}
\\ [2mm]
 & p^+ & \frac{m(2k+m)(m+k)}{2m+k} \left(\frac{3m+k}{m+k}\right)^\frac{k}{2k+m}  \left(\frac{2m}{m+k}\right)^\frac{m}{2k+m} \\ [2mm]
 & p^- & \frac{m(2k+m)(3m+k)}{2m+k} \left(\frac{m+k}{3m+k}\right)^\frac{k}{2k+m}  \left(\frac{2m}{3m+k}\right)^\frac{m}{2k+m}
\\[2mm] \hline\rule{0pt}{14pt}
\multirow{2}{*}{$W_3$} & q^\pm &
\frac{m(2k+m)(4m+2k+2-mq)}{2m+k+1} q^\frac{m}{2k+m}
\\ [2mm]
 &  p^\pm &
\frac{(p+1)m(2k+m)(3m^2+4km+4m+k^2+2k+1)}{p(2m+k+1)^2} p^{\frac{k}{2k+m}} \left(\frac{m (p+1)}{2m+k+1}\right)^{\frac{m}{2k+m}}
\\[2mm] \hline\rule{0pt}{14pt}
\multirow{4}{*}{$W_5$} & q^+ & \frac{l(2+l)}{4} 2^\frac{l-2}{l+2}\\ [2mm]
  &  q^- &
 \frac{(2+l)(l^2+2l-4)}{4l}\left(\frac{2(l-2)}{l}\right)^\frac{l-2}{l+2}
 \\ [2mm]
 l\geq 5 &  p^+ & \frac{l(2+l)}{4} \left(\frac{3l-4}{l}\right)^\frac{2}{l+2} \left(\frac{2(l-2)}{l}\right)^\frac{l-2}{l+2} \\ [2mm]
 &  p^- &  \frac{(3l-4)(2+l)}{4} \left(\frac{l}{3l-4}\right)^\frac{2}{l+2}  \left(\frac{2(l-2)}{3l-4}\right)^\frac{l-2}{l+2}
 \\[2mm] \hline \rule{0pt}{14pt}
\multirow{2}{*}{$W_6$} &  p^+ & \frac{28}{5}(120)^\frac{2}{7} \approx 21.9907\\ [2mm]
 & p^- & \frac{28}{5}(120)^\frac{2}{7} \approx 21.9907
\\[2mm] \hline \rule{0pt}{14pt}
\multirow{2}{*}{$W_{12}$} &  q^+ & \approx 69.1037 \\ [2mm]
 &  q^- & \approx 68.5187
\\[2mm] \hline \rule{0pt}{14pt}
\multirow{2}{*}{$W_{14}$} & p^+ & \approx 14.5750 \\ [2mm]
 & p^- & \approx 14.5750
\\[2mm] \hline \hline
\end{array}
$$}
\caption{Volume normalized scalar curvature of Einstein metrics in Table \ref{W3}.  In cases $W_1$ and $W_3$ we set $q:=q^\pm$ and $p:=p^\pm$.}
\label{W3Sc}
\end{table}

\begin{table}	
{\small 
$$
\begin{array}{c|c|c|c|c|c}
W_i & g & \lambda_\pg & \lambda_\pg^{max} & 2\rho &\text{Type}\\
[2mm]  \hline \hline
\rule{0pt}{14pt}
\multirow{2}{*}{$W_1$} &  g_{q^+} & \frac{5}{12} & \frac{5}{3} & \frac{2}{3} &  \text{saddle} \\ [2mm]
& g_{\times} &  &  &  &  \text{saddle}
\\[2mm] \hline \rule{0pt}{14pt}
\multirow{2}{*}{$W_1$}  & g_{q^+}  & \multirow{2}{*}{$\frac{2k+1}{k^2(k+1)}$} & \multirow{2}{*}{$\frac{2k+1}{k+1}$}  & \multirow{2}{*}{$\frac{k}{k+1}$} & \multirow{2}{*}{$\text{saddle}$}\\ [2mm] 
& k\geq 3, m=2 &&&& \\ [2mm]
 \hline\rule{0pt}{14pt}
\multirow{7}{*}{$W_1$}  & g_{q^+}  & \frac{m(4-q^2)}{2(k+2m-2)q} & \frac{q(2k+m)}{2(k+2m-2)}  & \frac{2(k+2m-2)-m q}{2(k+2m-2)} & \text{saddle} \\ [2mm]
  & \multirow{2}{*}{$g_{q^-}, T>0$}  & \multirow{2}{*}{$\frac{m(4-q^2)}{2(k+2m-2)q}$} & \multirow{2}{*}{$\frac{q(2k+m)}{2(k+2m-2)}$} & \multirow{2}{*}{$\frac{2(k+2m-2)-m q}{2(k+2m-2)}$} & \multirow{2}{*}{\text{{\bf loc.min.}}} \\ [2mm]
  & k<m &&&& \\ [2mm]
   & \multirow{2}{*}{$g_{q^-}, T>0$} & \multirow{2}{*}{$\frac{q(2k+m)}{2(k+2m-2)}$} & \multirow{2}{*}{$\frac{m(4-q^2)}{2(k+2m-2)q}$} & \multirow{2}{*}{$\frac{2(k+2m-2)-m q}{2(k+2m-2)}$} & \multirow{2}{*}{\text{{\bf loc.min.}}} \\ [2mm]
   & k>m &&&& \\ [2mm]
  & g_{q^-}, T<0 & \frac{q(2k+m)}{2(k+2m-2)} & \frac{m(4-q^2)}{2(k+2m-2)q} & \frac{2(k+2m-2)-m q}{2(k+2m-2)} & \text{saddle} \\ [2mm]
 &  g_{p^\pm} & {\text see \; \eqref{lamp}} & {\text see \; \eqref{lamp}} & \frac{(1+p^\pm)(1-4b^2)}{2p^\pm}& \text{saddle}\\ [2mm]
 \hline\rule{0pt}{14pt}
\multirow{5}{*}{$W_2$}& g_{q^+} & 0 & \frac{2k+m}{k+2m} & \frac{k+m}{k+2m} & \text{saddle}\\ [2mm]
 &  g_{q^-}, k<m & \frac{k}{k+m} & \frac{m(2k+m)}{(k+2m)(k+m)} & \frac{k^2+m^2+3km}{(k+m)(k+2m)} &  \text{{\bf loc.min.}} \\ [2mm]
 &  g_{q^-}, k>m&  \frac{m(2k+m)}{(k+2m)(k+m)}& \frac{k}{k+m} & \frac{k^2+m^2+3km}{(k+m)(k+2m)} &  \text{{\bf loc.min.}} \black\\ [2mm]
 & g_{p^+} & 0 & \frac{(k+5m)(k+m)}{(2m+k)(k+3m)} & \frac{m+k}{k+2m}& \text{saddle} \\ [2mm]
 & g_{p^-} & 0 & \frac{k+5m}{k+2m} & \frac{k+3m}{k+2m} & \text{saddle}
\\[2mm] \hline\rule{0pt}{14pt}
\multirow{4}{*}{$W_3$} & g_{q^+} &
 \frac{m(4-q^2)}{2q(k+2m+1)}& \frac{q(2k+m)}{2(k+2m+1)}& \frac{2(k+2m+1)-mq}{2(k+2m+1)} & \text{saddle} \\ [2mm]
 &  g_{q^-}, k<m & \frac{m(4-q^2)}{2q(k+2m+1)} & \frac{q(2k+m)}{2(k+2m+1)} & \frac{2(k+2m+1)-mq}{2(k+2m+1)}& \text{{\bf loc.min.}} \\ [2mm]
  &  g_{q^-}, k>m & \frac{q(2k+m)}{2(k+2m+1)} & \frac{m(4-q^2)}{2q(k+2m+1)}& \frac{2(k+2m+1)-mq}{2(k+2m+1)}& \text{{\bf loc.min.}} \\ [2mm]
 &  g_{p^\pm} & {\text see \; \eqref{lamp}} & {\text see \; \eqref{lamp}} & \frac{(1+p^\pm)(1-4b^2)}{2p^\pm}& \text{saddle}\\ [2mm] 
\hline\rule{0pt}{14pt}
\multirow{4}{*}{$W_5$} & g_{q^+} & 0 & \frac{l+2}{2(l-1)} & \frac{l}{2(l-1)}& \text{saddle}\\ [3mm]
  & g_{q^-}  & \frac{2}{l} & \frac{l^2-4}{2l(l-1)} &  \frac{l^2+2l-4}{2l(l-1)} &  \text{{\bf loc.min.}}\\ [3mm]
  & g_{p^+}  & 0 & \frac{l(5l-8)}{2(3l^2-7l+4)} & \frac{l}{2(l-1)} & \text{saddle} \\ [3mm]
  & g_{p^-}  & 0 & \frac{5l-8}{2(l-1)} & \frac{3l-4}{2(l-1)} & \text{saddle}
\\[2mm] \hline \hline
\end{array}
$$}
\caption{Einstein metrics in Table \ref{W3} on the spaces $W_1$-$W_5$ (see \eqref{Tkm} for the value of $T$ in case $W_1$).  They are all $G$-unstable and $G$-non-degenerate.  Any $g_{q^-}$ is a {\bf local minimum} with coindex $2$ and the rest are all saddle points with coindex $1$, except $g_{q^+}$ on $W_1$ ($m=2$), which has coindex $3$ if $k=2$ and coindex $2$ if $k\geq 3$ (see Examples \ref{Stiefel} and \ref{Stiefel3}).  Note that $g_{q^-}$ is the non-K\"ahler Einstein metric on the flag manifolds $W_2$ and $W_5$, the other three are K\"ahler-Einstein.}
\label{W4}
\end{table}

\begin{table}	
{\small 
$$
\begin{array}{c|c|c|c|c|c}
W_i & g & \lambda_\pg & \lambda_\pg^{max} & 2\rho &\text{Type}\\
[2mm]  \hline \hline
\rule{0pt}{14pt}
\multirow{2}{*}{$W_6$}  &  g_{p^+} &
\frac{67-\sqrt{1465}}{120} \approx 0.2393   &  \frac{67+\sqrt{1465}}{120}\approx 0.8772  & \frac{3}{5}=0.6&  \text{saddle}\\ [2mm]
 & g_{p^-} &\frac{67-\sqrt{1465}}{72}\approx 0.3989 &\frac{67+\sqrt{1465}}{72}\approx 1.4621 & 1 &  \text{saddle}
\\[2mm] \hline \rule{0pt}{14pt}
\multirow{2}{*}{$W_{12}$} &  g_{q^+} &
\frac{9-\sqrt{29}}{14}\approx 0.2582 & \frac{165+11\sqrt{29}}{210}\approx 1.0677 &\frac{55-\sqrt{29}}{70}\approx 0.7087 &  \text{saddle} \\ [2mm]
  & g_{q^-} & \frac{165-11\sqrt{29}}{210}\approx 0.5036 & \frac{9+\sqrt{29}}{14}\approx 1.0275 & \frac{55+\sqrt{29}}{70}\approx 0.8626 &  \text{saddle}
\\[2mm] \hline \rule{0pt}{14pt}
\multirow{2}{*}{$W_{14}$}& g_{p^+} & 0.1494  & 0.8657 & \frac{28(99+\sqrt{1177})}{9(499+\sqrt{1177})}\approx 0.5134 &  \text{saddle} \\ [2mm]
  & g_{p^-} & 0.3080  & 1.7839 & \frac{28(99-\sqrt{1177})}{9(499-\sqrt{1177})}\approx 1.0579&  \text{saddle}
\\[2mm] \hline \hline
\end{array}
$$}
\caption{Einstein metrics in Table \ref{W3} on the spaces $W_6$, $W_{12}$ and $W_{14}$.  They are all $G$-unstable and $G$-non-degenerate saddle points with coindex $1$.}
\label{W4-2}
\end{table}

For the family $W_1$ with $l=m$, it is easy to see that
\begin{equation}\label{Tkm}
T = \tfrac{-2k^2+2(k+m)(m-2)^2+8(m-1)}{(k+2m-2)^3}, 
\end{equation}
and
\begin{equation}\label{qpFam1}
p^{\pm}=\tfrac{5m^3+(9k-16)m^2+(20+5k^2-20k)m+k^3-6k^2-8+12k \pm 2(2m+k-2)\sqrt{D_1(k,m)}}{(k-2+m)(k-2+3m)(k+m)},
\end{equation}
where
$$
D_1(k,m):=\unm(k-1+m)(k+2m-2)^3T,
$$
which is positive if and only if $T>0$.   

On the other hand, for $W_3$ with $l=m$,
\begin{equation}\label{pmmFam3}
p^\pm=  \tfrac{10mk+5m^3+5mk^2+9m^2k+(k+1)^3+5m+8m^2 \pm (2m+k+1)\sqrt{D_2(k,m)}}{(k+1+3m)(k+1+m)(m+k)},
\end{equation}
where
$$
D_2(k,m):=  (2k+1+2m)((k+1)^2+4mk+2m^2k+2m^3+4m^2+4m),
$$
is always positive.  All these metrics have been listed  in Table \ref{W3} and the values of $\scalar_N(g)$ are given in Table \ref{W3Sc}.

We compute next the spectrum of $\lic_\pg$ for each of these metrics.  For $q=q^{\pm}$, we have that $2 \rho=1-bq$ by \eqref{rick}, and from \eqref{Lpgen} we obtain that
$$
L_\pg(g_{q}) = \left[\begin{matrix}
\frac{2b}{q} &  \tfrac{b(q^2-2)}{q}  & - \sqrt{bc} q \\
\tfrac{b(q^2-2)}{q} &   \frac{2b}{q} & - \sqrt{bc} q\\
-\sqrt{bc} q & - \sqrt{bc} q &  2cq
\end{matrix}\right], \qquad \Spec(\lic_\pg) = \left\{ 0, \tfrac{b(4-q^2)}{q}, q(b+2c)\right\}.
$$
For $p=p^{\pm}$, one obtains from \eqref{rick} that $2 \rho= \frac{(1+p)(1-4b^2)}{2p}$ and it follows from \eqref{Lpgen} that
$$
L_{\pg}(g_p) = \left[\begin{matrix}
\frac{p}{p+1}&  \frac{4b^2p^2+8b^2p+4b^2-p^2-1}{2p(p+1)}  & - \frac{\sqrt{c}(4b^2p+p+4b^2-1)}{2p\sqrt{b}} \\
\frac{4b^2p^2+8b^2p+4b^2-p^2-1}{2p(p+1)}  &  \frac{1}{p(p+1)}& -\frac{\sqrt{c}(4b^2p-p+4b^2+1)}{2p\sqrt{b}} \\
 - \frac{\sqrt{c}(4b^2p+p+4b^2-1)}{2p\sqrt{b}} & -\frac{\sqrt{c}(4b^2p-p+4b^2+1)}{2p\sqrt{b}} &  \frac{4bc(p+1)}{p}
\end{matrix}\right].  
$$
It is therefore straightforward to check that the spectrum is given by 
\begin{equation}\label{lamp}
\Spec(\lic_\pg) = \left\{ 0, \tfrac{4cb^2(p+1)^2 +b(p^2+1) \pm \sqrt{2bD(p)}}{2bp(p+1)}\right\},
\end{equation}
where $D(p):= c_4p^4 +c_3p^3 +c_2p^2 +c_3p+c_4$ and
$$
\begin{array}{c}
  c_4= (b+c)(8b^3c+1-4b^2+8b^4),\qquad c_3 = 8b^2(b+c)(4bc+4b^2-1), \\ [2mm]
  c_2= 96cb^4-8b^3+48c^2b^3+48b^5-2c-8cb^2.
\end{array}
$$
The case $W_1$ with $m=1$, i.e., the Stiefel manifold $\SO(k+2)/\SO(k)$, $k\geq 2$, was separately worked out in Example \ref{Stiefel} since it is not multiplicity-free.  

Finally, a straightforward inspection of all the above information gives that either $\lambda_\pg<2\rho<\lambda_\pg^{max}$ or $\lambda_\pg^{max}<2\rho$ in all cases.  This produces Table \ref{W4}, containing the $G$-stability and critical point types of all the metrics in Table \ref{W3}.     

\begin{remark}
For each space in Table \ref{W2}, the Einstein metrics $g_i$'s have the same value of $\scalar_N$ (see Table \ref{W2Sc}) and identical eigenvalues of $\lic_\pg$ up to scaling.   The same holds for the two Einstein metrics on $W_6$ and $W_{14}$ (see Table \ref{W3Sc}).  As far as we know, the question of whether these metrics are homothetic (i.e., isometric up to scaling) or not remains open, except for few particular cases.   
\end{remark}

\begin{table} {\small 	
$$
\begin{array}{c|c|c|c|c|c}
W_i &\ggo/\kg &  g & \lambda_\pg & \lambda_\pg^{max} & 2\rho \\
[2mm]  \hline \hline
\rule{0pt}{14pt}
\multirow{4}{*}{$W_2$}& \multirow{4}{*}{$\frac{\sug(k+l+m)}{\sg(\ug(k)\oplus\ug(l)\oplus\ug(m))}$}  & g_0 & \text{see\;\eqref{lam5}} & \text{see\;\eqref{lam5}} & \text{see\;\eqref{rho5}} \\[2mm]
&& g_k & 0 & \tfrac{4k+l+m}{(k+l+m)(2k+l+m)} & \frac{1}{k+l+m}
\\[2mm]
&& g_l & 0 & \tfrac{k+4l+m}{(k+l+m)(k+2l+m)} & \frac{1}{k+m+l}
\\[2mm]
&& g_m & 0 & \tfrac{k+l+4m}{(k+l+m)(k+l+2m)} & \frac{1}{k+m+l}
\\[2mm]
\hline\rule{0pt}{14pt}
\multirow{2}{*}{$W_8$}& \multirow{2}{*}{$\frac{\eg_6}{\spg(3)\oplus\spg(1)}$} & g_1=(1, 1.4618, 1.8845) & 0.1605 &  0.9669 & 0.5745 \\[2mm]
 && g_2=(1, 0.8640, 0.4838) & 0.3464 & 1.6227 & 1.0116
\\[2mm] \hline \rule{0pt}{14pt}
\multirow{2}{*}{$W_{10}$}&\multirow{2}{*}{$\frac{\eg_7}{\sug(6)\oplus \spg(1) \oplus \RR}$} & g_1=(1, 0.8882, 0.5717) & 0.4354 & 1.3150 & 0.9492\\[2mm]
&&g_2=(1, 1.1896, 1.6291) & 0.2118  & 1.0217 & 0.6480
\\[2mm] \hline \hline
\end{array}
$$}
\caption{Einstein metrics on some generalized Wallach spaces with $a_i\ne a_j$ for $i\ne j$ (see  \cite{ChnKngLng} and \cite{FrsLmsNkn}).  In case $W_2$, $k,l,m$ are pairwise different and for $W_8$ and $W_{10}$ the numbers are all approximate values.   All these metrics are $G$-unstable and $G$-non-degenerate.  The non-K\"ahler metric $g_0$ on the flag manifold $W_2$ (note that the other three are K\"ahler-Einstein) has coindex $2$ (i.e., $\lambda_\pg^{max}<2\rho$), so it is a {\bf local minimum}, and the rest are all saddle points of coindex $1$ (i.e., $\lambda_\pg<2\rho<\lambda_\pg^{max}$).  See Example \ref{W4exa} for $W_4$.}
\label{W5}
\end{table}

\begin{table}
{\small	
$$
\begin{array}{c|c|c|c}
G/K & g & (x_1,x_2,x_3) & (\det_{g_{\kil}}g)^{\frac{1}{n}}\scalar(g)  \\
[2mm]  \hline \hline
\rule{0pt}{14pt}
\multirow{3}{*}{$\frac{E_8}{E_6 \times \SU(2)\times \U(1)}$} & g_0 & (1,2,3) &   \frac{1577}{30}1207959552^\frac{1}{83} \approx 67.6289 \\ [2mm]
& g_1 & (1, 0.914286, 1.54198) & 66.9159 \\ [2mm]
& g_2 & (1, 1.0049, 0.129681) & 65.6151
\\[2mm] \hline\rule{0pt}{14pt}
\multirow{3}{*}{$\frac{E_8}{\SU(8) \times \U(1)}$} & g_0 & (1,2,3) & \frac{782}{15} 1152^\frac{1}{23} \approx 70.8307 \\ [2mm]
& g_1 & (1, 0.717586, 1.25432) & 69.5453 \\ [2mm]
& g_2 & (1, 1.06853, 0.473177) & 69.1155
\\[2mm] \hline\rule{0pt}{14pt}
\multirow{3}{*}{$\frac{E_7}{\SU(5)\times \SU(3) \times \U(1)}$} & g_0 & (1,2,3) & \frac{250}{9} {24}^\frac{1}{10} \approx 38.1696 \\ [2mm]
& g_1 & (1, 0.678535, 1.201221) & 37.4141 \\ [2mm]
& g_2 & (1, 1.090568, 0.546044) & 37.3277
\\[2mm] \hline\rule{0pt}{14pt}
\multirow{3}{*}{$\frac{E_7}{\SU(6)\times \SU(2) \times \U(1)}$} & g_0 & (1,2,3) & \frac{517}{18} {294912}^\frac{1}{47} \approx 37.5486 \\ [2mm]
& g_1 & (1, 0.85368, 1.45259) & 37.0717 \\ [2mm]
& g_2 & (1, 1.01573, 0.229231) & 36.4084
\\[2mm] \hline\rule{0pt}{14pt}
\multirow{3}{*}{$\frac{E_6}{\SU(3)\times \SU(3) \times \SU(2) \times \U(1)}$} & g_0 & (1,2,3) & \frac{203}{12} {4608}^\frac{1}{29} \approx 22.6278  \\ [2mm]
& g_1 & (1, 0.771752, 1.33186) & 22.2677 \\ [2mm]
& g_2 & (1, 1.04268, 0.373467) & 22.0134
\\[2mm] \hline\rule{0pt}{14pt}
\multirow{3}{*}{$\frac{F_4}{\SU(3)\times \SU(2) \times \U(1)}$} & g_0 & (1,2,3) & \frac{100}{9} 24^\frac{1}{10} \approx 15.2678  \\ [2mm]
& g_1 & (1, 0.678535, 1.201221) & 14.9656 \\ [2mm]
& g_2 & (1, 1.090568, 0.546044) & 14.9311
\\[2mm] \hline\rule{0pt}{14pt}
\multirow{3}{*}{$\frac{G_2}{\SU(2)\times\U(1)}$} & g_0 & (1,2,3) & \frac{25}{12} 18^\frac{1}{5} \approx 3.7137 \\ [2mm]
& g_1 & (1, 1.67467, 2.05238) & 3.7104 \\ [2mm]
& g_2 & (1, 0.186894, 0.981478) & 3.4422
\\[2mm] \hline \hline
\end{array}
$$}
\caption{Einstein metrics on flag manifolds with $b_2(M)=1$ and $r=3$ (see \cite{AnsChr2}).  On the fourth row, the embedding of $\SU(6)\times\SU(2)\times\U(1)$ on $E_7$ is different from the one for the generalized Wallach space $W_{10}$.}
\label{FM3}
\end{table}

\begin{table}
{\small
$$
\begin{array}{c|c|c|c}
G/K & g & (x_1,x_2,x_3,x_4) & (\det_{g_{\kil}}g)^{\frac{1}{n}}\scalar(g) \\
[2mm]  \hline \hline
\rule{0pt}{14pt}
\multirow{3}{*}{$\frac{F_4}{\SU(3)\times \SU(2) \times \U(1)}$} & g_0 & (1,2,3,4) & \frac{70}{9} 24^\frac{1}{10} \approx 14.5995  \\ [2mm]
& g_1 & (1,1.2761, 1.9578, 2.3178) & 14.5693 \\ [2mm]
& g_2 & (1,0.9704, 0.2291, 1.0097) & 14.0370
\\[2mm] \hline\rule{0pt}{14pt}
\multirow{3}{*}{$\frac{E_7}{\SU(4)\times \SU(3)\times \SU(2)\times \U(1)}$} & g_0 &
(1,2,3,4) &  \frac{212}{9} 24^\frac{8}{53} \approx 38.0563 \\ [2mm]
& g_1 & (1, 0.8233, 1.2942, 1.3449) & 37.6284 \\ [2mm]
& g_2 & (1, 0.9912, 0.5783, 1.1312) & 37.3618
\\[2mm] \hline\rule{0pt}{14pt}
\multirow{3}{*}{$\frac{E_{8}}{\SU(7)\times \SU(2)\times \U(1)}$} & g_0 & (1,2,3,4) & \frac{637}{15}\sqrt{2} 3^\frac{1}{7} \approx 70.2624\\ [2mm]
& g_1 & (1, 0.9133, 1.4136, 1.5196) & 69.6567\\ [2mm]
& g_2 & (1, 0.9663, 0.4898, 1.0809) & 68.8049
\\[2mm] \hline\rule{0pt}{14pt}
\multirow{5}{*}{$\frac{E_{8}}{\SU(10)\times \SU(3)\times \U(1) }$} & g_0 & (1,2,3,4) & \frac{658}{15} 41472^\frac{2}{47} \approx 68.9659 \\ [2mm]
& g_1 & (1,0.6496, 1.1094, 1.0610) & 66.5752 \\ [2mm]
& g_2 & (1, 1.1560, 1.0178, 0.2146) & 66.1753 \\ [2mm]
& g_3 & (1, 1.0970, 0.7703, 1.2969) & 66.9855 \\ [2mm]
& g_4 & ( 1, 0.7633, 1.0090, 0.1910) & 65.7898
\\[2mm] \hline \hline
\end{array}
$$}
\caption{Einstein metrics on flag manifolds with $b_2(M)=1$ and $r=4$ (see \cite{ArvChr}).}
\label{FM4}
\end{table}

\begin{table}	
{\small
$$
\begin{array}{c|c|c|c}
G/K & g & (x_1, \dots, x_r) & (\det_{g_{\kil}}g)^{\frac{1}{n}}\scalar(g) \\
[2mm]  \hline \hline
\rule{0pt}{14pt}
\multirow{6}{*}{$\frac{E_8}{\SU(5)\times \SU(4)\times\U(1)}$} & \multirow{2}{*}{$g_0$} & \multirow{2}{*}{$(1,2,3,4,5)$} & \multirow{2}{*}{$\frac{572}{15} 2^\frac{25}{52} 3^\frac{5}{26} 5^\frac{1}{26}$} \\ [2mm] 
&&& \approx 69.9307  \\ [2mm]
& g_1 & (1, 0.599785, 1.08371, 0.901823, 1.22291) & 68.8905 \\ [2mm]
& g_2 & (1, 1.02137, 0.546007, 1.05352, 1.10879) & 68.7023  \\ [2mm]
& g_3 & (1, 1.08294, 1.04088, 0.532615, 1.10351) & 68.7757 \\ [2mm]
& g_4 & (1, 0.720713, 1.02546, 0.475234, 1.07095) & 68.6913 \\ [2mm]
& g_5 & (1, 1.03732, 1.04718, 1.03082, 0.29862) & 68.4798
\\[2mm] \hline\rule{0pt}{14pt}
\multirow{5}{*}{$\frac{E_8}{\SU(5)\times \SU(3) \times \SU(2)\times\U(1)}$} & \multirow{2}{*}{$g_0$} & \multirow{2}{*}{$(1,2,3,4,5,6)$} & \multirow{2}{*}{$\frac{159}{5} 2^\frac{65}{106} 3^\frac{25}{106} 5^\frac{5}{53}$}\\ [2mm]  
&&&\approx 69.0420\\ [2mm]
& g_1 & (1, 0.823084, 1.1467, 1.17377, 1.42664, 1.46519) & 68.6856 \\ [2mm]
& g_2 & (1, 0.986536, 0.636844, 1.06853, 1.13323, 0.921127) & 68.4684 \\ [2mm]
& g_3 & (1, 0.90422, 0.778283, 0.927483, 1.03408, 0.359949) & 68.2283 \\ [2mm]
& g_4 & (1, 0.954875, 0.965321, 1.00534, 0.290091, 1.01965) & 67.8054
\\[2mm] \hline \hline
\end{array}
$$}
\caption{Einstein metrics on flag manifolds with $b_2(M)=1$ and $r=5,6$ (see \cite{ChrSkn1}).}
\label{FM56}
\end{table}

\begin{table}	
{\small
$$
\begin{array}{c|c|c|c|c}
G/K & g &\lambda_\pg & \lambda_\pg^{max} &   2\rho  \\
[2mm]  \hline \hline
\rule{0pt}{14pt}
\multirow{3}{*}{$\frac{E_8}{E_6 \times \SU(2)\times \U(1)}$} & g_0 & \frac{6}{5}-\frac{1}{30}\sqrt{51} & \frac{6}{5}+\frac{1}{30}\sqrt{51} & \frac{19}{30}    \\ [2mm]
& g_1 & 0.478572 & 1.55965 & 0.821452  \\ [2mm]
& g_2 & 0.118731 & 1.272251 & 0.829109
\\[2mm] \hline\rule{0pt}{14pt}
\multirow{3}{*}{$\frac{E_8}{\SU(8) \times \U(1)}$} & g_0  & \frac{11 -\sqrt{6}}{10} & \frac{11 +\sqrt{6}}{10} & \frac{17}{30}  \\ [2mm]
& g_1  & 0.4676547 & 1.340237 & 0.819949  \\ [2mm]
& g_2 & 0.326773 & 1.166896 & 0.785751
\\[2mm] \hline\rule{0pt}{14pt}
\multirow{3}{*}{$\frac{E_7}{\SU(5)\times \SU(3) \times \U(1)}$} & g_0  & \frac{5}{6} & \frac{4}{3} & \frac{5}{9}  \\ [2mm]
& g_1 & 0.469601 & 1.311724 & 0.825338 \\ [2mm]
& g_2 & 0.352429 & 1.157063 & 0.772758
\\[2mm] \hline\rule{0pt}{14pt}
\multirow{3}{*}{$\frac{E_7}{\SU(6)\times \SU(2) \times \U(1)}$} & g_0  &  \frac{7-\sqrt{2}}{6} & \frac{7+\sqrt{2}}{6}& \frac{11}{18} \\ [2mm]
& g_1 & 0.469565 & 1.485343 & 0.816530  \\ [2mm]
& g_2 & 0.194860 & 1.232011 & 0.820660
\\[2mm] \hline\rule{0pt}{14pt}
\multirow{3}{*}{$\frac{E_6}{\SU(3)\times \SU(3) \times \SU(2) \times \U(1)}$} & g_0 & \frac{9}{8}-\frac{1}{24}\sqrt{33}& \frac{9}{8}+\frac{1}{24}\sqrt{33}& \frac{7}{12}   \\ [2mm]
& g_1 & 0.465849  & 1.391997 & 0.815861  \\ [2mm]
& g_2 & 0.281933 & 1.187358 & 0.801967
\\[2mm] \hline\rule{0pt}{14pt}
\multirow{3}{*}{$\frac{F_4}{\SU(3)\times \SU(2) \times \U(1)}$} & g_0 & \frac{5}{6} & \frac{4}{3} & \frac{5}{9}  \\ [2mm]
& g_1 & 0.469601 & 1.311722 & 0.825338 \\ [2mm]
& g_2 & 0.352428 & 1.157064 & 0.772757
\\[2mm] \hline\rule{0pt}{14pt}
\multirow{3}{*}{$\frac{G_2}{\U(2)}$} & g_0 & \unm & \frac{5}{4} & \frac{5}{12} \\ [2mm]
& g_1 & 0.413430 & 1.211009 & 0.502068  \\ [2mm]
& g_2 & 0.19355 & 2.670881 & 0.970058
\\[2mm] \hline \hline
\end{array}
$$}
\caption{Einstein metrics in Table \ref{FM3}.  The K\"ahler metric $g_0$ is always {\bf $G$-stable} and the rest are all $G$-non-degenerate and $G$-unstable of coindex $1$.}
\label{FS3}
\end{table}

\begin{table}
{\small 	
$$
\begin{array}{c|c|c|c|c|c|c}
G/K & g &\lambda_\pg & \lambda_2 & \lambda_\pg^{max} &  2\rho & \text{cx} \\
[2mm]  \hline \hline
\rule{0pt}{14pt}
\multirow{3}{*}{$\frac{F_4}{\SU(3)\times \SU(2) \times \U(1)}$} & g_0  &  \frac{1}{2} & \frac{20}{27} & \frac{7}{6}  & \frac{14}{36}&0 \\ [2mm]
& g_1  & 0.4111 & 0.8447 & 1.2064 & 0.5380 &1 \\ [2mm]
& g_2  & 0.2108 & 1.2583 & 2.1506 & 0.8231 &1
\\[2mm] \hline\rule{0pt}{14pt}
\multirow{3}{*}{$\frac{E_7}{\SU(4)\times \SU(3)\times \SU(2)\times \U(1)}$} & g_0 &
  \frac{11-\sqrt{13}}{12} & \frac{53}{54} & \frac{11 +\sqrt{13}}{12} & \frac{4}{9} &0 \\ [2mm]
& g_1  & 0.5001 & 1.0457 & 1.2838 & 0.7173 &1\\ [2mm]
& g_2 & 0.4331 & 1.0696 & 1.3332 & 0.7626 &1
\\[2mm] \hline\rule{0pt}{14pt}
\multirow{3}{*}{$\frac{E_{8}}{\SU(7)\times \SU(2)\times \U(1)}$} & g_0 & \frac{9}{10}-\frac{\sqrt{21}}{15} & \frac{14}{15} & \frac{9}{10}+\frac{\sqrt{21}}{15}
& \frac{13}{30} & 0 \\ [2mm]
& g_1  & 0.4826 & 1.0134 & 1.2468 & 0.6781 &1\\ [2mm]
& g_2  & 0.3937  & 1.1307 & 1.3923 & 0.7826 &1
\\[2mm] \hline\rule{0pt}{14pt}
\multirow{5}{*}{$\frac{E_{8}}{\SU(10)\times \SU(3)\times \U(1) }$} & g_0 & \frac{19}{20}-\frac{\sqrt{309}}{60} & \frac{97}{90} & \frac{19}{20}+\frac{\sqrt{309}}{60} & \frac{7}{15} & 0\\ [2mm]
& g_1  & 0.5039 & 1.1437  & 1.4399  & 0.7970 &1\\ [2mm]
& g_2  & 0.1698 & 0.9025  & 1.1115 & 0.7038  &1 \\ [2mm]
& g_3  & 0.4527 & 1.0347 & 1.3018 & 0.7173 &1\\ [2mm]
& g_4  & 0.2157  & 0.6048 & 1.5116 & 0.8030 & {\bf 2}
\\[2mm] \hline \hline
\end{array}
$$}
\caption{Einstein metrics in Table \ref{FM4}.  The K\"ahler metric $g_0$ is always {\bf $G$-stable} and all the other ones are $G$-non-degenerate and $G$-unstable with coindex as in the last column.}
\label{FS4}
\end{table}

\begin{table}	
{\small
$$
\begin{array}{c|c|c|c|c|c|c}
 g & \lambda_\pg & \lambda_2 & \lambda_3 &\lambda_\pg^{max} & 2 \rho & cx \\
[2mm]  \hline \hline
\rule{0pt}{14pt}
 g_0 & 0.483308 & 0.807779 & 1.002382 & 1.156529 &  \frac{11}{30}\approx 0.36  &0 \\ [2mm]
 g_1 & 0.510214 & 1.027352 & 1.238460 & 1.655126 & 0.757540  &1 \\ [2mm]
 g_2 & 0.444857 & 1.112644 & 1.251940 & 1.598781 & 0.731014  &1 \\ [2mm]
 g_3 & 0.387659 & 0.850534 & 1.094764 & 1.233562 & 0.678788  &1\\ [2mm]
 g_4 & 0.458617 & 0.613228 & 1.437954 & 1.612703 & 0.773963  &{\bf 2}\\ [2mm]
 g_5 & 0.240674 & 0.791323 & 1.169060 & 1.320873 & 0.674542  & 1
\\[2mm] \hline \hline
\end{array}
$$}
\caption{Einstein metrics on $E_8/ \SU(5)\times\SU(4)\times\U(1)$ (see Table \ref{FM56}).  The K\"ahler metric $g_0$ is {\bf $G$-stable} and the remaining ones are all $G$-non-degenerate and $G$-unstable.}
\label{FS5}
\end{table}

\begin{table}
{\small 	
$$
\begin{array}{c|c|c|c|c|c|c|c}
 g & \lambda_\pg & \lambda_2 & \lambda_3 & \lambda_4 &\lambda_\pg^{max} & 2 \rho & cx \\
[2mm]  \hline \hline
\rule{0pt}{14pt}
 g_0 & 0.373056 & 0.616209 & 0.784713 & 0.907323 & 1.025363 & \frac{3}{10}=0.3 & 0 \\ [2mm]
 g_1 & 0.560023 & 0.753518 & 0.990611 & 1.177858 & 1.218291 & 0.627866 & 1 \\ [2mm]
 g_2 & 0.499591 & 0.995988 & 1.060532 & 1.188804 & 1.256011 & 0.697205 & 1 \\ [2mm]
 g_3 & 0.371562 & 0.622626 & 0.981804 & 1.317906 & 1.460064 & 0.735036 &  {\bf 2} \\ [2mm]
 g_4 & 0.262279 & 0.832143 & 1.186281 & 1.385267 & 1.473554 & 0.698590 & 1
\\[2mm] \hline \hline
\end{array}
$$}
\caption{Einstein metrics on $E_8/ \SU(5)\times \SU(3)\times \SU(2)\times\U(1)$ (see Table \ref{FM56}).  The K\"ahler metric $g_0$ is {\bf $G$-stable} and the other ones are all $G$-unstable.}
\label{FS6}
\end{table}

\subsection{Case $d_1,d_2,d_3$ pairwise different}\label{dist2-sec}
This case will be treated partially, since the Einstein metrics are not available in the literature for all the spaces.  It follows from \cite{ChnKngLng} that in the case $W_2$ (with $k,l,m$ pairwise different), $G/K$ admits exactly four $G$-invariant Einstein metrics up to scaling (see also \cite[Theorem 5]{FrsLmsNkn}):
$$
\begin{array}{lcl}
g_0=(l+m,k+m,k+l),&& g_k=(l+m+2k,k+m,k+l), \\ [2mm]
g_l=(l+m,k+2l+m,k+l),  && g_m=(l+m,k+m,k+l+2m).
\end{array}
$$
For $g_0$ one obtains from \eqref{rick} that
\begin{equation}\label{rho5}
2\rho= \tfrac{(k+l)(k+m)(l+m)+2klm}{(k+l+m)(k+l)(k+m)(l+m)},
\end{equation}
and from \eqref{Lpgen} that
\begin{equation}\label{lam5}
\Spec(\lic_\pg(g_0))=\left\{ 0, \tfrac{(k+l)(k+m)(l+m) + 4klm \pm \sqrt{P(k,l,m)}}{2(k+l+m)(k+l)(k+m)(l+m)}\right\},
\end{equation}
where
$$
P(k,l,m) = (k+l)(k+m)(l+m)\left((k+l)(k+m)(l+m)-8klm\right).
$$
For the others metrics on $W_2$,
as well as for the cases $W_8$ and $W_{10}$, which both admit exactly two Einstein metrics (see \cite{ChnKngLng}), we proceeded as above and all the information has been collected in Table \ref{W5}.  The volume normalized scalar curvature of the metrics on $W_8$ and $W_{10}$ are given by 
\begin{align*}
\scalar_N(W_8,g_1) \approx 21.7434, \qquad  \scalar_N(W_8,g_2) \approx 21.5470, \\
\scalar_N(W_{10},g_1) \approx 36.7796, \qquad  \scalar_N(W_{10},g_2) \approx 37.1468.
\end{align*}

\begin{example}\label{W4exa} 
It is proved in \cite[Section 4.1]{ChnKngLng} that $W_4=\SU(2l)/\U(l)$, $l\geq 2$, admits exactly two Einstein metrics, given by $g=(1,x_2,x_3)$, where $x_3$ is one of the only two positive roots of the quartic polynomial 
$$
12 l^4x^4
- (48 l - 8)l^3x^3 + (72 l^2 - 36 l - 4)l^2x^2
- (48 l^3 - 48l^2 + 4 l + 4)lx + 12 l^4 - 20 l^3 + 7l^2 + 2l - 1,
$$
and 
$$
x_2= \tfrac{2l^2x_3^2+2lx_3+ 1-l-2l^2}{2l(2lx_3-2l+1)}.
$$
Using respectively \eqref{Lpgen} and \eqref{rick}, it is straightforward to check that both metrics satisfy: 
$$
\{\lambda_\pg, \lambda_\pg^{max}\} = \tfrac{(l-1)x_2^2 + (l+1)x_3^2 + l \pm \sqrt{Q(l,x_2,x_3)}}{4lx_2 x_3}, \qquad 
2 \rho=\tfrac{1+4x_2 x_3-x_2^2-x_3^2}{4x_2x_3}, 
$$
where 
$$
Q(l,x_2,x_3)= (2l+4l^2)x_3^4+(2-4l^2-4l^2x_2^2+2l)x_3^2
-1+2x_2^2-2lx_2^2-4l^2x_2^2+4l^2x2^4-2lx_2^4+4l^2.
$$
With the help of a computer, we have checked that $\lambda_\pg < 2\rho < \lambda_\pg^{max}$ holds for any $2\leq l\leq 20$, so in these cases both metrics are $G$-non-degenerate and $G$-unstable saddle critical points of $\scalar|_{\mca_1^G}$ of coindex $1$.  This numerical evidence suggests that these metrics are saddle points for any $l\geq 2$.   
\end{example}

The only cases which have not been considered in this paper all have pairwise different $d_1,d_2,d_3$ and correspond to the spaces $W_1$ (between two and four Einstein metrics, see \cite[Theorem 1.3]{ChnNkn}) and $W_3$ (four Einstein metrics, see \cite[Example 4]{FrsLmsNkn}).  This is due to the fact that the Einstein metrics are not explicitly given in these papers, only existence results are provided.

\subsection{Application to the prescribed Ricci curvature problem}\label{RLI}
Given a symmetric $2$-tensor $T$ on a differentiable manifold $M$, asking about the existence and uniqueness (up to scaling) of a Riemannian metric $g$ and a constant $c>0$ such that
$
\ricci(g) = cT,
$
is called the prescribed Ricci curvature problem (see e.g.\ \cite[Chapter 5]{Bss}).  For a $G$-invariant tensor $T$ and metric $g$ on a homogeneous space $M=G/K$, the question is therefore about the image and the injectivity (up to scaling) of the differentiable function
$
\ricci:\mca^G\longrightarrow\sca^2(M)^G.
$
A metric $g\in\mca^G$ is called {\it Ricci locally invertible} if $\ricci$ is as invertible as it can be near $g$, i.e., for any $T$ near $\ricci(g)$, existence and (local) uniqueness of solutions to the prescribed Ricci curvature problem are guaranteed (see \cite{PRP} for a more detailed treatment).   According to \cite[Lemma 6.1]{PRP} (see also \cite[Lemma 4.5]{stab-tres}), $d\ricci|_g = \unm\lic_\pg$ after identifications, so a necessary condition for the Ricci local invertibility of $g$ is that $\dim{\Ker\lic_\pg}|_{\sca^2(M)^G}=1$ (see \cite[Theorem 1.3]{PRP}).  The metrics with $\lambda_\pg=0$ contained in Tables \ref{W2}, \ref{W4} and \ref{W5} are therefore examples of $G$-invariant Einstein metrics which are not Ricci locally invertible.  These are precisely the K\"ahler-Einstein metrics on flag manifolds with $b_2(M)=2$.  On the other hand, it is easy to check that the remaining Einstein metrics considered in this paper are all Ricci locally invertible.  

In what follows, we give examples of curves of $G$-invariant metrics with identical Ricci tensors.  The metric $g_{KE}=(1,1,2)$ is K\"ahler-Einstein on each of the following homogeneous spaces,
$$
\SO(2l)/\U(1)\times\U(l-1), \quad l\geq 4, \qquad E_6/\SO(8)\times \U(1)\times \U(1), 
$$
and $g_{KE}=(l+m,k+m,k+l+2m)$ is a K\"ahler-Einstein metric on
$$
\SU(k+l+m)/\Se(\U(k)\times\U(l)\times\U(m)), \qquad k,l,m\geq 1.  
$$
Note that the last $g_{KE}$ is a multiple of $(1,1,2)$ if $k=l$.  In all the above cases, a straightforward computation using \eqref{rick} gives that if $g=(x_1,x_2,x_3)$, then $\ricci(g)=\ricci(g_{KE})$ if and only if $x_3=x_1+x_2$ (recall that the Ricci tensor in terms of a $-\kil_\ggo$-orthonormal basis is given by $\ricci(g)=(x_1\rho_1,x_2\rho_2,x_3\rho_3)$).  This implies that the metrics 
$$
g_t:=\left(t, \unm(-t+(t^2+\tfrac{8}{t})^{1/2}), \unm(t+(t^2+\tfrac{8}{t})^{1/2})\right), \quad t>0,
$$   
have all the same volume and Ricci tensor as $g_1=g_{KE}=(1,1,2)$.  This curve was recently discovered by Pulemotov and Ziller on the full flag manifold $\SU(3)/\Se(\U(1)^3)$.  

On the spaces $\SO(8)/\U(1)\times\U(3)$, $E_6/\SO(8)\times \U(1)\times \U(1)$ and $\SU(3k)/\Se(\U(k)^3)$ (i.e., when $d_1=d_2=d_3$), the scalar curvature is given by 
$$
\scalar(g_t) = c\left(-t^2+3(t^2+\tfrac{8}{t})^{1/2}+\tfrac{4}{t}\right),\qquad\mbox{for some}\; c>0,  
$$
which is a convex function on $(0,\infty)$ with a global minimum at $t=1$.  The family $\{ g_t:1\leq t\}$ is therefore pairwise non-homothetic as $\scalar(g_t)$ is strictly increasing on $[1,\infty)$ (note that also $\{ g_t:0<t\leq 1\}$ is pairwise non-homothetic).  Since for any $t>1$, $g_t$ has the same Ricci eigenvalues as $g_s$, where $s:=\unm\left(-t+(t^2+\tfrac{8}{t})^{1/2}\right)$, we do not know whether these pairs of metrics are isometric or not.   

For the rest of the K\"ahler-Einstein metrics which are not Ricci locally invertible, i.e., $g_{p^{\pm}}$ on $\SO(2l)/\U(1)\times\U(l-1)$, $l\geq 5$ and $W_2$ with $k=l$, and the three K\"ahler-Einstein metrics on $W_2$ with pairwise different $k,l,m$, the same behavior occurs and an explicit curve $g_t$ such that $\ricci(g_t)\equiv\ricci(g_1)$ and with the same properties as above can easily be given.  

\begin{remark}\label{new2}
After the first version of this paper appeared in arXiv, these examples were generalized in \cite{PlmZll}.  
\end{remark}

\section{Flag manifolds}\label{flag-sec}

A {\it flag manifold} (or generalized flag manifold) is a homogeneous space $M=G/K$, where $G$ is a compact semisimple Lie group and $K$ is the centralizer of a torus in $G$.  It is called {\it full flag} when $K$ is a maximal torus of $G$.  Any compact, simply connected and de
Rham irreducible homogeneous K\"ahler manifold is isometric to a $G$-invariant metric on a flag manifold $G/K$ with $G$ simple and simply connected.  The isotropy representation of any flag manifold is multiplicity-free.  We refer to \cite{Arv,AnsChr} and references therein for more information on flag manifolds.    

In this section, we compute the $G$-stability and critical point types of all Einstein metrics on flag manifolds with $b_2(M)=1$, where $b_2(M)$ denotes the second Betti number of $M$.  As in \S\ref{W-sec}, this will follow from the knowledge of the spectrum of $\lic_\pg$, which is computed using Theorem \ref{formLp} and explicitly provided for all Einstein metrics.  Our results on $G$-stability agree with \cite[Theorem 3.1]{AnsChr}, where the Ricci flow behavior on flag manifolds with $b_2(M)=1$ is studied, and also with \cite[Theorems 1 and 2]{Grm}, where the case of flag manifolds with three isotropy summands is considered.  

The number $r$ of $\Ad(K)$-irreducible summands of the isotropy representation of these spaces goes from $2$ to $6$, and each of them admits a unique K\"ahler-Einstein metric given by $g_0:=(1,2,\dots,r)$, with respect to a suitable decomposition $\pg=\pg_1\oplus\dots\oplus\pg_r$ in $\Ad(K)$-irreducible subspaces.  Note that $\dim{\mca_1^G}=r-1$.  

The case when $r=2$ has been worked out in Example \ref{r2-flag}.  Actually any flag manifold with $r=2$ has $b_2(M)=1$.  On the other hand, the remaining flag manifolds with $r=3$ all have $b_2(M)=2$ and are given by the generalized Wallach spaces $W_2$, $W_5$ and $W_7$, already studied in \S\ref{W-sec}.  

The structural constants and numerical approximations of all Einstein metrics (other than $g_0$) on flag manifolds with $b_2(M)=1$ are given in \cite{AnsChr2} for $r=3$ (the case $E_7/\SU(5)\times\SU(3)\times\U(1)$ had to be corrected since the right dimensions are $(60,30,10)$, see \cite{ChrSkn2}), in \cite{ArvChr} for $r=4$ and in \cite{ChrSkn1} for $r=5,6$.  The metrics, including their volume normalized scalar curvature $\scalar_N$ relative to the Killing metric $g_{\kil}=(1,\dots,1)$ (see \eqref{scalN}), are listed in Tables \ref{FM3}, \ref{FM4} and \ref{FM56}.  It follows that the Einstein metrics on each space are always pairwise non-homothetic.  Also note that $\scalar_N(g_0)>\scalar_N(g_i)$ for any $i\geq 1$ in all cases.  

It is straightforward to check that the corresponding numerical approximations for the eigenvalues of the $r\times r$ matrix $\lic_\pg$ and the $G$-stability types of these Einstein metrics are given as in Tables \ref{FS3}, \ref{FS4}, \ref{FS5} and \ref{FS6}, according to the number $r$.  They are all $G$-non-degenerate and Ricci locally invertible (see \S\ref{RLI}).


\section{Appendix: Maple codes}

In this appendix, we include the computational Maple codes used in the paper.   

The first part concerns generalized Wallach spaces and computes the three eigenvalues (one of them is always zero) of the generic $3\times 3$ matrix of $\lic_\pg$ given in \eqref{Lpgen}, in terms of the structural constants $a_i:=[123]/d_i$, $i=1,2,3$ (see Table \ref{W1}) and the Einstein metric $g=(x_1,x_2,x_3)$.  The eigenvalues exhibited in Tables \ref{W2}, \ref{W4}, \ref{W4-2} and \ref{W5} are obtained by just replacing in the general formula for each case with the corresponding $a_i$'s and $x_i$'s, which are given in Tables \ref{W2}, \ref{W3} and \ref{W5}.     

In the second part, we consider flag manifolds and first compute the Ricci eigenvalues and the entries of the generic matrix of $\lic_\pg$ in the case $r=6$, i.e., for $E_8/ \SU(5)\times \SU(3)\times \SU(2)\times\U(1)$.  It turns out that for the other cases $r=3,4,5$ such a matrix is obtained by considering the principal $r\times r$ submatrix and setting $Aijk=0$ whenever at least one of $i,j,k$ is $>r$.   

As an example, the third part computes all the information given in the last row of Table \ref{FM56} and Table \ref{FS6} for the five Einstein metrics on $E_8/ \SU(5)\times \SU(3)\times \SU(2)\times\U(1)$.


\begin{figure}
   \includegraphics[width=150mm]{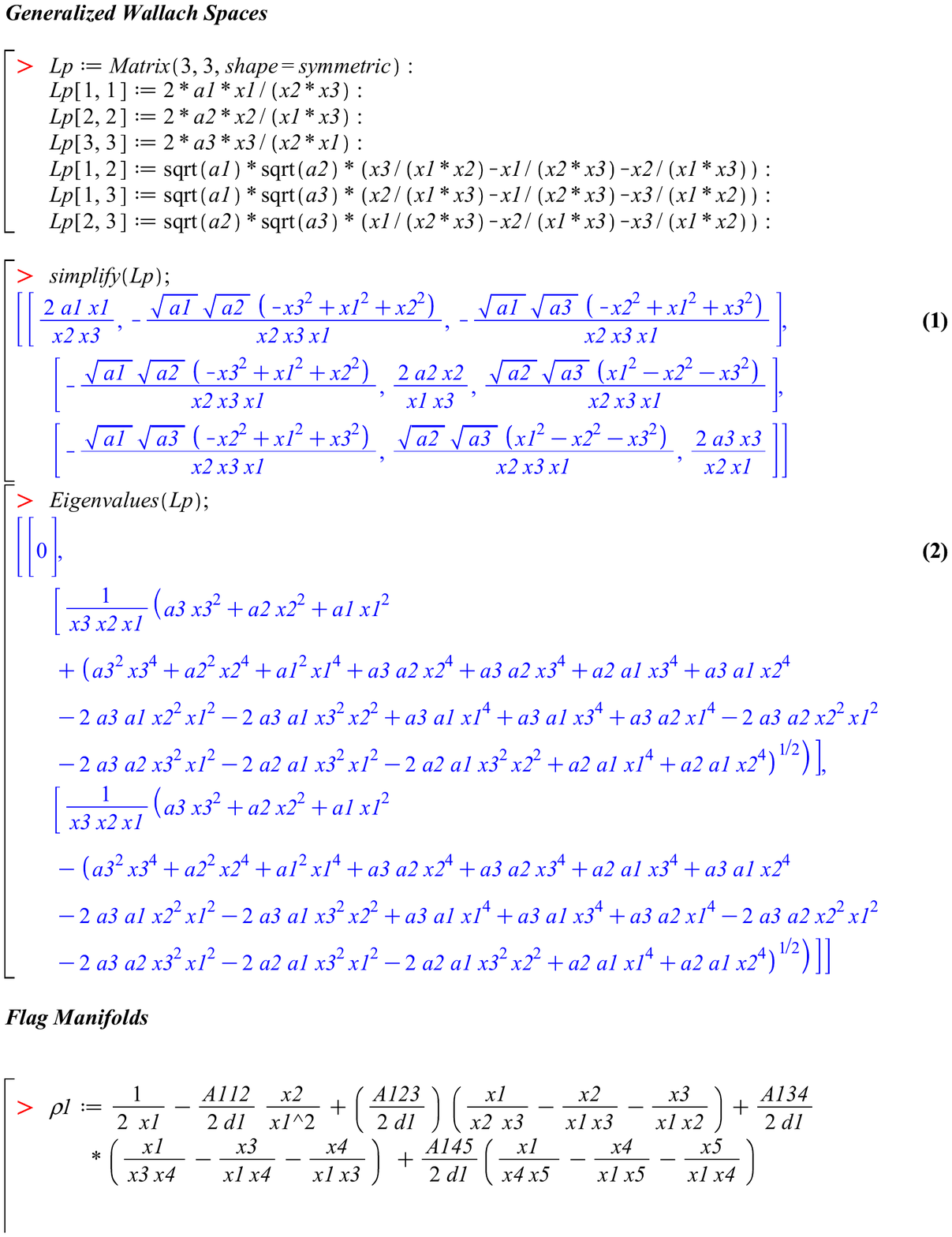}
    \end{figure}

\begin{figure}
   \includegraphics[width=150mm]{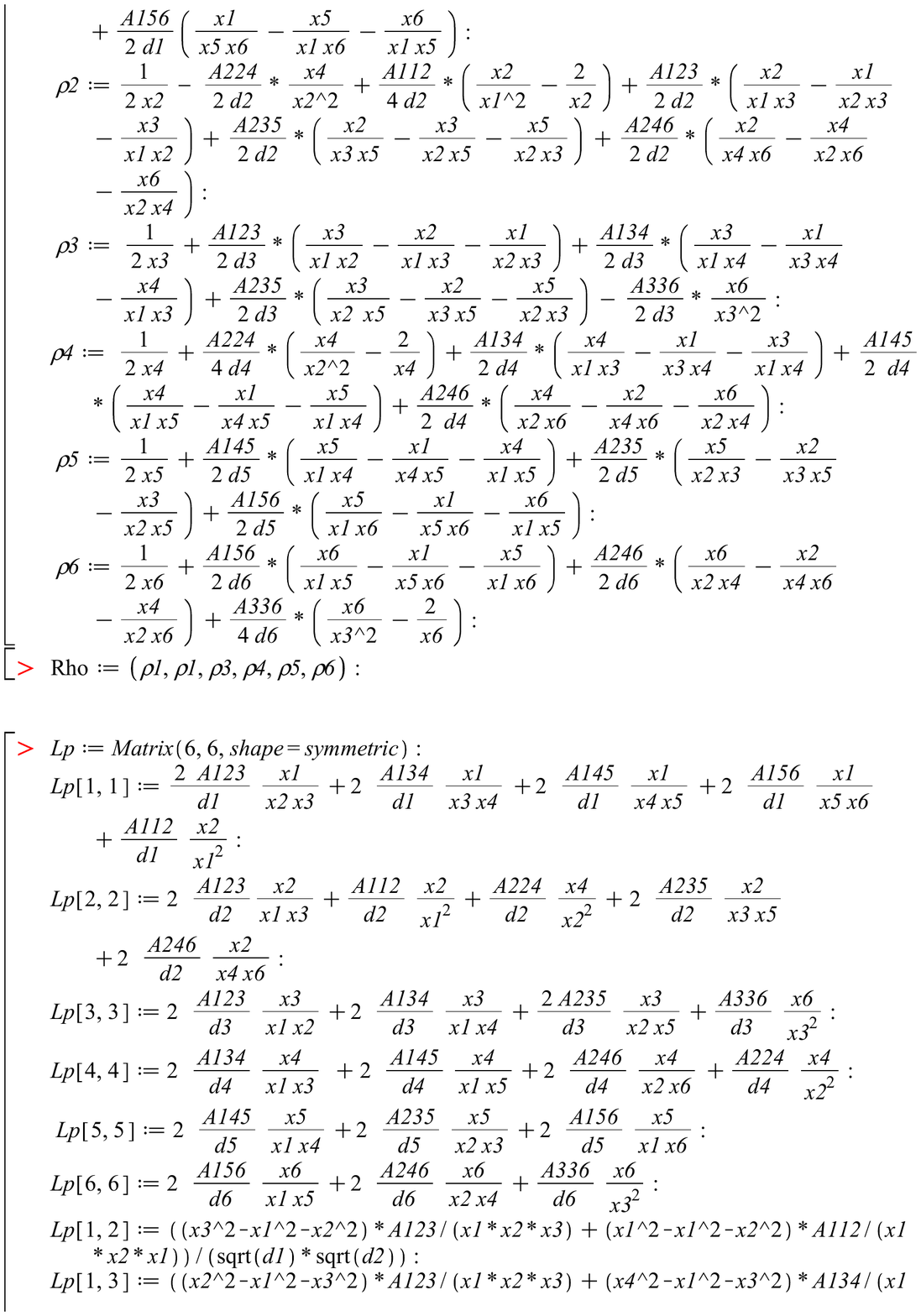}
    \end{figure}
    
\begin{figure}
   \includegraphics[width=150mm]{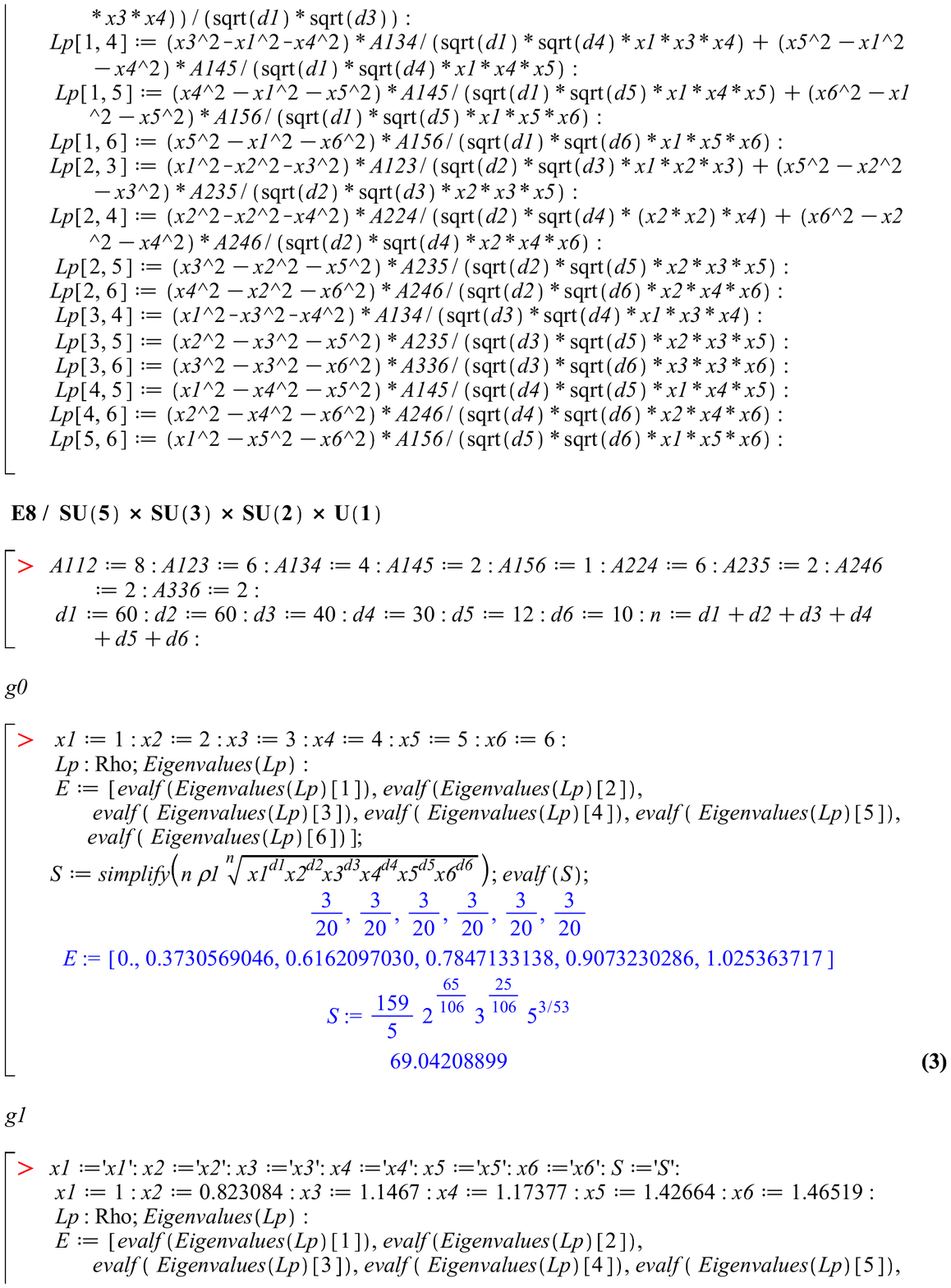}
    \end{figure}                

\begin{figure}
   \includegraphics[width=150mm]{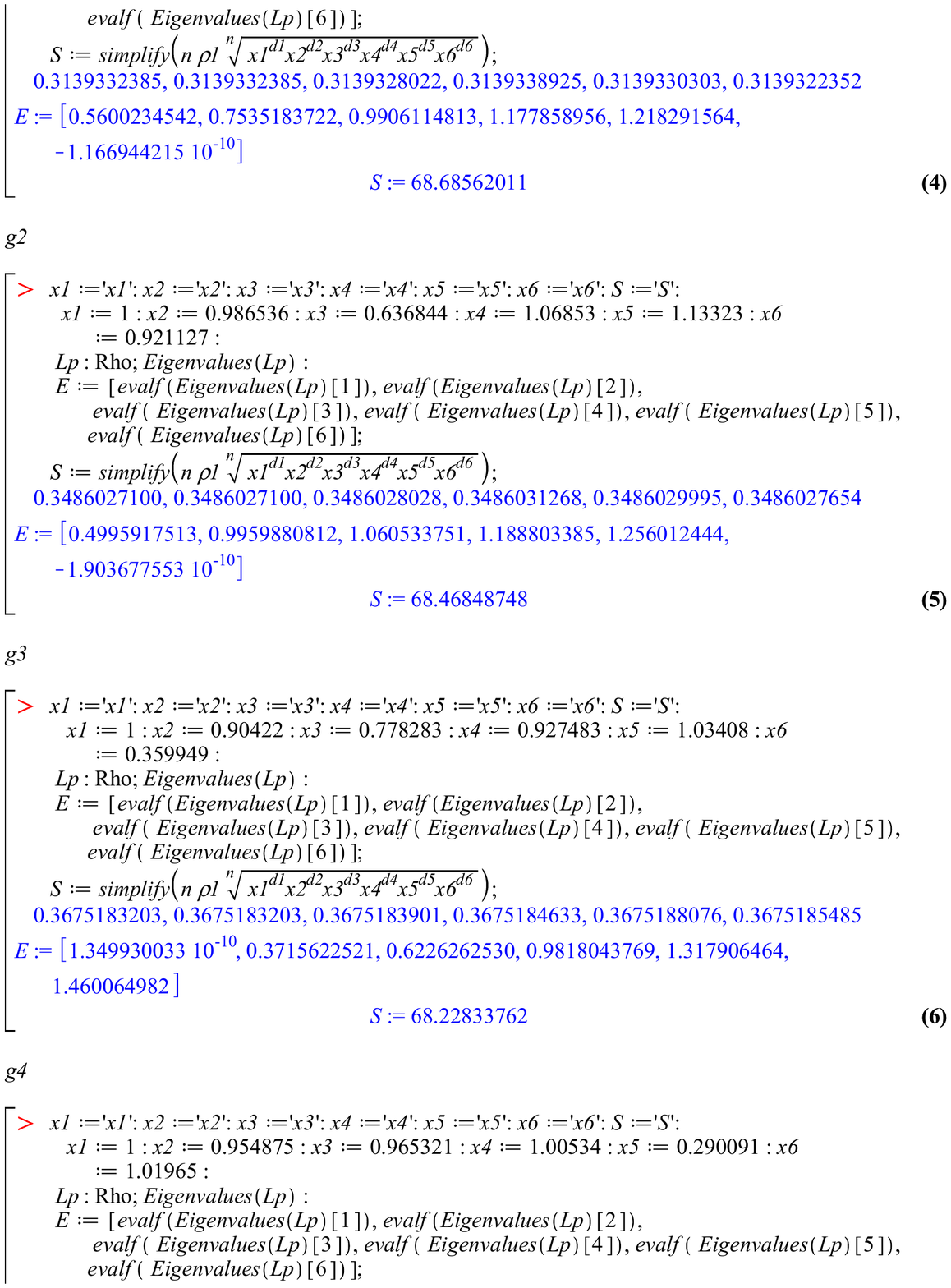}
    \end{figure}
    
\begin{figure}
   \includegraphics[width=150mm]{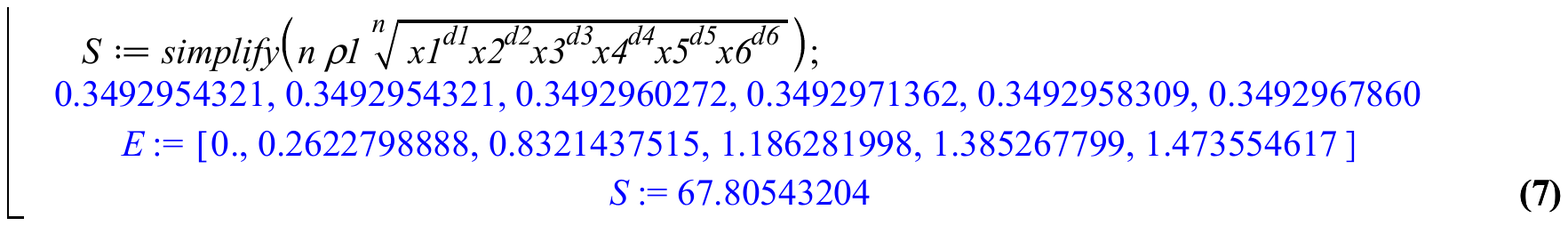}
    \end{figure}

\end{document}